\documentclass[12pt]{article}

\usepackage{amsmath}
\usepackage{amsfonts}
\usepackage{amssymb}
\usepackage{amsthm}
\usepackage{graphicx}

\usepackage{cite}
\usepackage{caption}
\usepackage{url}
\usepackage{color}\usepackage{placeins}
\usepackage{svg}

\usepackage{algorithm}
\usepackage{algpseudocode}
\usepackage{mathtools}
\usepackage{subcaption}
\usepackage{multirow}
\usepackage{bbm}
\usepackage{epstopdf}
\usepackage{graphicx}
\usepackage{bm}
\usepackage{enumerate}
\usepackage{booktabs}

\theoremstyle{definition}

\newtheorem{remark}{Remark}
\newcommand{\norm}[1]{\left\lVert#1\right\rVert}

\algdef{SE}{Begin}{End}{\textbf{begin}}{\textbf{end}}

\title{General-domain FC-based shock-dynamics solver I:\\ Basic elements}

\author{Oscar P. Bruno\footnote{Computing and Mathematical Sciences, Caltech,
    Pasadena, CA 91125, USA} \and Daniel V. Leibovici\footnote{NVIDIA Research, Santa Clara, CA, 95051, USA}}
    
\begin{document}
\date{}
\maketitle


\begin{abstract}
  This contribution, Part~I in a two-part article series, presents a
  general-domain version of the {\em FC-SDNN} (Fourier Continuation
  Shock-detecting Neural Network) spectral scheme for the numerical
  solution of nonlinear conservation laws, which is applicable under
  arbitrary boundary conditions and in {\em general domains}.  Like
  the previous simple-domain contribution (Journal of Computational
  Physics X {\bf 15}, (2022)), the present approach relies on the use
  of the Fourier Continuation (FC) method for accurate spectral
  representation of non-periodic functions (without the limiting CFL
  time-step constraints inherent in other spectral schemes) in
  conjunction with smooth artificial viscosity assignments localized
  in regions detected by means of a Shock-Detecting Neural Network
  (SDNN). Relying on such techniques, the present Part~I paper
  introduces a novel multi-patch/subpatch artificial viscosity-capable
  domain decomposition strategy for {\em complex domains with smooth
    boundaries}, and it illustrates the methodology by means of a
  variety of computational results produced by an associated parallel
  implementation of the resulting shock-capturing algorithm in a
  present-day computing cluster. The subsequent Part~II contribution
  then extends the algorithm to enable treatment of obstacles with
  non-smooth boundaries, it considers questions concerning
  parallelization and accuracy, and it presents comparisons with
  physical theory and prior experimental and computational results.
  The resulting multi-patch FC-SDNN algorithm does not require use of
  problem-dependent algorithmic parameters or positivity-preserving
  limiters, and, on account of its use of an overlapping-patch
  discretization, it is geometrically flexible and efficiently
  parallelized. A variety of numerical tests for the 2D Euler
  equations are presented, including the simulation of supersonic and
  hypersonic flows and shocks past physical obstacles at high speeds,
  such as e.g. Mach 25 re-entry flow speeds.
\end{abstract}

\vspace{0.5 cm}
\maketitle
\newpage
\section{\label{sec:introduction}Introduction}

In a sequence of two papers (whose Part~I is the present contribution)
we present a general-domain {\em FC-SDNN} (Fourier Continuation
Shock-detecting Neural Network) spectral scheme for the numerical
solution of nonlinear conservation laws under arbitrary boundary
conditions and in {\em general domains}. Like its simple-domain
counterpart~\cite{bruno2022fc}, the approach proposed here enjoys
limited spatial dispersion away from shocks and other discontinuities
and does not suffer from the limiting CFL constraints inherent in
other spectral schemes. The present Part~I paper introduces a novel
multi-patch/subpatch artificial viscosity-capable domain decomposition
strategy for {\em complex domains with smooth boundaries}, and it
illustrates the methodology by means of a variety of computational
results produced by an associated parallel implementation of the
resulting shock-capturing algorithm in a present-day computing
cluster. The Part~II contribution~\cite{bruno_leibo_partII}, in turn,
extends the algorithm to enable treatment of obstacles with non-smooth
boundaries; it includes results validating the accuracy and energy
conservation properties of the proposed multi-patch-based artificial
viscosity method; and it introduces a parallel implementation of high
parallel efficiency and scalability.

The proposed approach relies on use of the FC-Gram Fourier
Continuation method~\cite{albin2011spectral,
  amlani2016fc,bruno2010high, lyon2010high} for spectral
representation of non-periodic functions in conjunction with new
multi-patch versions of the Shock-Detecting Neural Network (SDNN) and
smooth artificial viscosity methods introduced
in~\cite{schwander2021controlling,bruno2022fc}. The neural network
approach, which was introduced in~\cite{schwander2021controlling},
utilizes Fourier series to discretize the gas dynamics and related
equations in two-dimensional periodic rectangular domains, eliminates
Gibbs ringing at shocks (whose locations are determined by means of an
artificial neural network) by assigning artificial viscosity to a
small number of discrete points in the immediate vicinity of the
shocks.  The use of the classical Fourier spectral method in that
contribution restricts the method's applicability to periodic problems
(thus precluding, in particular, applicability to geometries
containing physical boundaries), and its highly localized viscosity
assignments give rise to a degree of non-smoothness, resulting in
certain types of unphysical oscillations manifested as serrated
contour-lines in the flow fields.  As demonstrated
in~\cite{bruno2022fc}, the FC-SDNN approach leverages smooth viscosity
assignments to deliver accurate, essentially oscillation-free
solutions for non-periodic problems with arbitrary boundary
conditions. Building on this foundation, the present two-part
contribution extends the methodology to general geometries via
multi-domain decompositions.  Illustrative results include Euler
simulations of re-entry-speed flows past obstacles---reaching, for
instance, Mach 25---highlighted in
Section~\ref{subsubsec:flow_cylinder_inf}; cf.
Remark~\ref{rem_hyperson}.

The computational solution of systems of conservation laws has been
tackled by means of a variety of numerical methods. The brief account
provided in~\cite[Sec. 1]{bruno2022fc} reviews some of the important
methodologies developed in this field over the last few decades,
including: Low-order finite volume and finite difference methods,
possibly equipped with slope limiters~\cite{leveque1992numerical,
  leveque2002finite}; ENO~\cite{harten1987uniformly} and
WENO~\cite{liu1994weighted, jiang1996efficient} schemes as well as a
hybrid FC/WENO solver~\cite{shahbazi2011multi}; Artificial
viscosity-based methods in the context of finite difference
discretizations~\cite{richtmyer1948proposed, vonneumann1950method,
  gentry1966eulerian,lapidus1967detached, lax1959systems} and
Discontinuous Galerkin discretizations~\cite{persson2006sub}; The
entropy-viscosity (EV) method~\cite{guermond2011entropy,
  kornelus2018flux}; The spectral viscosity
method~\cite{maday1989analysis,tadmor1989convergence,tadmor1990shock,tadmor1993super,maday1993legendre},
that relies on convolution with a compactly supported kernel and
which, while demonstrated numerically on the basis of 1D problems,
may, like the EV method, potentially result in oversmearing; The
C-method~\cite{reisner2013space,ramani2019space1, ramani2019space2},
which augments the hyperbolic system with an additional equation used
to determine a spatio-temporally smooth viscous term, and which, like
the EV method, requires selection of several problem-dependent
parameters and algorithmic variations; and, finally, a number of
ML-based
techniques~\cite{ray2018artificial,discacciati2020controlling,stevens2020enhancement,schwander2021controlling,
  shahane2024rational} which enhance the performance of classical
shock capturing schemes. The relative advantages provided by the
proposed FC-SDNN, which were reviewed in~\cite[Sec. 1]{bruno2022fc},
are also present in the general-domain multi-patch algorithms
introduced in the present Part~I-II paper sequence---resulting, in
particular, in accurate solutions with smooth contour lines for flows
in general domains.

As indicated above, the FC-SDNN method presented in this paper
provides a general solver for shock-wave problems governed by
conservation laws by incorporating 1)~The shock-detection methodology
mentioned above; in conjunction with, 2)~The aforementioned rapidly
convergent FC-Gram Fourier Continuation representations (which,
eliminating the Gibbs phenomenon that would otherwise arise from lack
of periodicity, provide rapidly convergent Fourier expansions of
non-periodic functions); as well as, 3)~Smooth artificial viscosity
assignments~\cite{bruno2022fc} based on use of certain newly-designed
{\em smooth viscosity windows}---which results in (a)~Smooth flow
profiles away from discontinuities, as well as (b)~Sharp shock
resolution, while (c)~Enabling applicability to gas-dynamics problems
including very strong shocks, without recourse to
positivity-preserving limiters for the density and pressure
fields. Unlike polynomial spectral methods such as the one based on
use of Chebyshev representations, the proposed approach does not give
rise to quadratically refined meshes near interval endpoints, and it
utilizes viscosity assignments that decay proportionally to the
spatial mesh size, resulting in, 4)~CFL restrictions that decay
linearly with the spatial meshsize. Finally, in view of its use of an
overlapping-patch domain decomposition
strategy~\cite{albin2011spectral,amlani2016fc,chesshire1990composite},
5)~The approach enjoys significant geometrical flexibility, afforded
by the loose connectivities required between adjacent
overlapping-patch computational subdomains, leading to a lesser
geometry-processing overhead than other discretization methods.

The FC-SDNN shock-capturing capability, which does not use
problem-dependent parameters, employs an Artificial Neural Network
(ANN) for detection of discontinuities together with smooth localized
viscosity assignments. Briefly, the approach relies on a single
pre-trained neural network for detection of discontinuities on the
basis of Gibbs oscillations in Fourier series, together with the
spectral resolution provided by the FC method in the fully general
context of non-periodic problems with given boundary
conditions---allowing for sharp resolution of flow features such as
shocks and contact discontinuities and smooth accurate evaluation of
flow fields away from discontinuities. In view of its innovative
smooth viscosity assignments, further, this procedure effectively
eliminates Gibbs oscillations while avoiding introduction of a type of
flow-field roughness that is often evidenced by the serrated contour
levels produced by other methods. And, importantly, as indicated in
Remark~\ref{adiabatic_th}, the algorithm does not require use of
positivity-preserving limiters to prevent the occurrence of negative
densities and pressures---even for flows past obstacles at large Mach
number values.

This paper is organized as follows. Section~\ref{Preliminaries}
presents necessary preliminaries concerning the hyperbolic problems
under consideration together with a brief introduction to the Fourier
Continuation method and the basic elements of the artificial-viscosity
strategies used. Section~\ref{sec:Multidomain} describes the proposed
overlapping-patch domain-decomposition and discretization strategy for
general complex domains. Sections~\ref{subsec:time_marching}
and~\ref{sec:viscosity} then describe algorithmic elements that enable
the implementation of multi-patch versions of the FC-SDNN
time-marching scheme and the artificial viscosity assignment method,
respectively; a brief outline of the overall algorithm is provided in
Section~\ref{sec:algorithm}. Section~\ref{sec:numerical results}
presents a number of applications illustrating the method's
performance, including flows and shocks past obstacles at hypersonic
speeds, all of which were obtained on the basis of a Message Passing
Interface (MPI) massive parallel implementation of the FC-SDNN
method---the details and properties of which are presented in
Part~II. A few concluding remarks, finally, are presented in
Section~\ref{sec:conclusion}.

\section{\label{Preliminaries} Preliminaries}

\subsection{\label{equations} Governing equations}

This paper proposes novel Fourier spectral methodologies for the
numerical solution of the two-dimensional Euler equations in general
non-periodic domains and with general boundary conditions. In detail,
we consider the 2D Euler equations
\begin{equation} \label{eq: euler 2d equation}
\dfrac{\partial}{\partial t}
\begin{pmatrix}
\ \rho \\[\jot]
\ \rho u\\[\jot]
\ \rho v\\[\jot]
\ E
\end{pmatrix} + \frac{\partial}{\partial x}
\begin{pmatrix}
\ \rho u\\[\jot]
\ \rho u^2 + p\\[\jot]
\ \rho u v\\[\jot]
\ u (E + p)
\end{pmatrix}
+ \frac{\partial}{\partial y}
\begin{pmatrix}
\ \rho v\\[\jot]
\ \rho u v\\[\jot]
\ \rho v^2 + p\\[\jot]
\ v (E + p)
\end{pmatrix} = 0
\end{equation}
on a  domain $\Omega\subset \mathbb{R}^2 $, where $\rho$,
\textit{\textbf{u}}, and $p$ denote the density, velocity
vector and  pressure, respectively, and where $E$ denotes
the total energy:
\begin{equation}\label{perfect_gas}
    E = \frac{p}{\gamma - 1} + \frac{1}{2} \rho \lvert \textit{\textbf{u}} \rvert ^2.
\end{equation}
Letting $\theta$ denote the temperature and noting that
\begin{equation}\label{p_T}
    p = \rho \theta,
\end{equation}
we may also write
\begin{equation}\label{perfect_gas_T}
  E = \frac{\rho \theta}{\gamma - 1} + \frac{1}{2} \rho \lvert \textit{\textbf{u}} \rvert ^2,
\end{equation}
where $\gamma$ denotes the heat capacity ratio; the ideal diatomic-gas
heat-capacity value $\gamma = 1.4$ is used in all of the test cases
considered in this paper.  The temperature $\theta$, which is not an
unknown in equation~\eqref{eq: euler 2d equation}, is mentioned here
in view of the important role it plays in regard to enforcement of
adiabatic boundary conditions---as indicated in
Remark~\ref{adiabatic}.

Letting $\mathbb{I}$ denote the identity tensor, denoting by
$\textbf{a}\otimes \textbf{b} = (a_ib_j)$ the tensor product of the
vectors $\textbf{a} = (a^i)$ and $\textbf{b} = (b^j)$, and calling
  \begin{equation}\label{flow_flux}
    \textbf{e} = (\rho, \rho \textbf{u}, E)^T\quad \mbox{and}\quad 
    f(\textbf{e}) = (\rho\textbf{u}, \rho
    \textbf{u} \otimes \textbf{u} + p\mathbb{I},
    \textbf{u}(E + p))^T,
\end{equation}
equation~\eqref{eq: euler 2d equation} may be expressed in the form
\begin{equation} \label{eq: nonlinear-eqn} \dfrac{\partial}{\partial
    t} \textbf{e}(\textbf{x}, t) + \nabla \cdot \big(
  f(\textbf{e}(\textbf{x}, t))\big) = 0,\quad \textbf{e} : \Omega \times [0, T] \rightarrow \mathbb{R}^4,
\end{equation}
where $f : \mathbb{R}^4 \rightarrow \mathbb{R}^4 \times \mathbb{R}^2$
denotes a (smooth) convective flux.  (The term
$\nabla \cdot \big(f(\textbf{e})\big)$ in~\eqref{eq:
  nonlinear-eqn} is a three coordinate vector whose first, second and
third coordinates are a scalar, a vector and a scalar,
respectively. Using the Einstein summation convention, these three
components are respectively given by
$\nabla\cdot(\rho\textbf{u}) = \partial_i (\rho u^i)$,
$\big( \nabla\cdot (\rho \textit{\textbf{u}} \otimes
\textit{\textbf{u}} + p\mathbb{I}) \big)^j =\partial_i (\rho u^j u^i +
p)$ and $\nabla\cdot((E + p)\textbf{u}) = \partial_i ((E + p) u^i)$.)

The speed of sound
\begin{equation}\label{snd_sp}
  a = \sqrt{\gamma \frac{p}{\rho}}
\end{equation}
for the Euler equations plays important roles in the various
artificial viscosity assignments considered in this paper.

\subsection{\label{fc} FC-based spatial approximation and differentiation}

Relying on the FC method, reference~\cite{bruno2022fc} presented a
low-dispersion solver based on the use of Fourier expansions for
equation~\eqref{eq: euler 2d equation} in non-periodic contexts. The
present paper further extends the applicability of the FC approach to
enable the use of mesh discretizations which, based on the use of
multiple patches and subpatches, is applicable to arbitrary spatial
domains. In the framework proposed in this paper, at each time-stage
of a multi-stage Runge-Kutta scheme
(Section~\ref{subsec:differentiation}) the equation is evolved
independently on each subpatch (Section~\ref{subsec:overl_dec}) in
curvilinear form (Section~\ref{subsec:curvilinear_eqns})---thus
effectively requiring the approximation of derivatives of functions
given on 2D Cartesian grids. Clearly this can be accomplished on the
basis of sequential applications of a 1D FC-based differentiation
procedure---such as the FC-Gram algorithm~\cite{albin2011spectral}
employed in this paper.

The FC-Gram algorithm constructs an accurate Fourier approximation of
a given generally non-periodic function $F$ defined on a given
one-dimensional interval---which, without loss of generality, is
assumed in this section to equal the unit interval $[0,1]$:
$F:[0,1]\to \mathbb{R}$.  As indicated in~\cite{bruno2022fc} and
references therein, and as briefly outlined in what follows, starting
from given function values $F_j =F(x_j)$ at $N$ points
$x_j=jh\in [0,1]$ ($h=1/(N-1)$), the FC-Gram algorithm produces a
function
\begin{equation}\label{FC_exp}
 F^c(x) = \sum^{M}_{k = -M} \hat{F}^c_k \exp (2 \pi i k x /\beta)
\end{equation}
which is defined (and periodic) in an interval $[0,\beta]$ that
strictly contains $[0,1]$, where $\hat{F}^c_k$ denote the FC
coefficients of $F$ and where, as detailed below, for $N$ large we
have $M\sim N/2$. To do this the FC-Gram algorithm first uses the
function values $(F_0,\dots,F_{d-1})^T$ and $(F_{N-d},...,F_{N-1})^T$
(where $d$ is a suitably small positive integer, typically satisfying
$2\leq d\leq 5$), that is, the values taken by $F$ in small matching
subintervals $[0,\Delta]$ and $[1-\Delta,1]$ of length
$\Delta = (d-1)h$ near the left and right endpoints of the interval
$[0,1]$, to produce, at first, a discrete (but ``smooth'') periodic
continuation vector $\mathbf{F}^c$ of the vector
$\mathbf{F} = (F_0,\dots ,F_{N-1})^T$. (Throughout this paper
the value $d=2$ has been used.)  The continuation $\mathbf{F}^c$
contains the entries in the vector $\mathbf{F}$ followed by a number
$C$ of continuation function values in the interval $[1,\beta]$, so
that the extension $\mathbf{F}^c$ transitions smoothly from $F_{N-1}$
back to $F_0$. (The FC method can also be applied to certain
combinations of function values and derivatives by constructing the
continuation vector $\mathbf{F}^c$ on the basis of e.g. the vector
$\mathbf{F} = (F_0,\dots,F_{N-2},F'_{N-1})^T$, where
$F_j\approx F(x_j)$ for $1\leq j\leq N-2$ and where
$F'_{N-1} \approx F'(x_{N-1})$.  Such a procedure enables imposition
of Neumann boundary conditions in the context of the FC method;
see~\cite[Sec. 6.3]{amlani2016fc}).) The resulting vector
$\mathbf{F}^c$ can be viewed as a discrete set of values of a smooth
and periodic function which can be used to produce the Fourier
continuation function $F^c$ via an application of the FFT algorithm.

As detailed in~\cite{bruno2022fc} and references therein, the
necessary discrete continuation values of the function $F^c$ at the
points $N/(N-1),(N+1)/(N-1),\ldots,(N+C-1)/(N-1)$ can be expressed in
the matrix form
\begin{equation}\label{eq: continuation}
    \mathbf{F}^c = 
    \begin{pmatrix}
        \mathbf{F} \\[\jot]
        A_{\ell} Q^T \mathbf{F}_{\ell} + A_r Q^T\mathbf{F}_r 
    \end{pmatrix}
\end{equation}
where the $d$-dimensional vectors $\mathbf{F}_{\ell}$ and
$\mathbf{F}_r$ contain the point values of $F$ at the first and last
$d$ discretization points in the interval $[0,1]$, respectively; where
$Q$ is certain $d \times d$ orthogonal matrix; and where $A_{\ell}$
and $A_r$ are certain (small) $C \times d$, matrices which can be
inexpensively computed once and stored on disc, and then read for use
to produce FC expansions for functions $G:[a,b]\to \mathbb{R}$ defined
on a given 1D interval $[a,b]$, via re-scaling to the interval
$[0,1]$.  Spatial derivatives, in turn, can be computed by
differentiation of FC-based Fourier expansions constructed by means of
the FC-Gram procedure, followed by evaluation at all necessary grid
points via application of the IFFT algorithm to the Fourier
coefficients
\begin{equation}\label{fc_der}
  (\hat{F}^{c})'_k= \frac{2 \pi i k}{\beta}\hat{F}^c_k.
\end{equation}

The order of accuracy of approximations obtained via the Fourier
Continuation method (whether derivative values or function values are
prescribed at endpoints) is restricted by the corresponding order $d$
of the Gram polynomial expansion, which was selected as $d=2$ for all
of the examples presented in this paper.  The relatively low order of
accuracy afforded by this selection, rendered necessary for stability
purposes in the context of interactions between shocks and physical
boundaries, is not a matter of consequence as regards to the problems
considered in the present paper, where high orders of accuracy are not
expected from any numerical method on account of shocks and other flow
discontinuities.

\subsection{\label{artvisc} Artificial Viscosity}

The presence of shocks and other flow discontinuities in the solutions
of equation~(\ref{eq: nonlinear-eqn}) and other nonlinear conservation
laws amounts to a significant computational challenge, usually
manifesting itself in spurious ``Gibbs oscillations'', low accuracies
and, even, incorrect numerical flow predictions. As is known
(cf.~\cite{schwander2021controlling,bruno2022fc} and references
therein), a carefully selected artificial viscosity term may be used
to induce solution smoothness while maintaining an adequate level of
numerical accuracy, particularly away from discontinuities. The latter
goal is best achieved by introducing viscous terms that vanish (or
nearly vanish) away from the solution discontinuities. In particular,
the contribution~\cite{bruno2022fc} utilizes a viscous term that is
not only localized, but also {\em smooth}, and which, thus, avoids
introduction of discretization level non-smoothness that otherwise
manifests itself, for example, in the existence of serrated level-set
curves.

The FC-SDNN framework~\cite{bruno2022fc}, which is extended here to a
multi patch-based general-domain discretization scheme (in which
viscosity smoothness is preserved across discretization patches, as
detailed in Section~\ref{sec:viscosity}, and more specifically in
Section~\ref{subsec:window_propagation}), utilizes a viscous term of
the form
\begin{equation}\label{f_visc}
f_\textit{visc}[\textbf{e}] = \mu[\textbf{e}]  \nabla(
  \textbf{e}(\textit{\textbf{x}}, t))
\end{equation}
given in terms of a certain ``viscosity operator''
$\mu[\textbf{e}](x,t)$, which is added to the right hand side of
(\ref{eq: nonlinear-eqn}), resulting in the viscous equation
\begin{equation} \label{eq: convection diffusion eqn}
  \frac{\partial \textbf{e}(\textit{\textbf{x}}, t)}{\partial t} + \nabla \cdot
  \big( f(\textbf{e}(\textit{\textbf{x}}, t))\big) = \nabla \cdot
  \big( f_\textit{visc}[\textbf{e}](\textit{\textbf{x}}, t)\big).
\end{equation}
The FC-SDNN method additionally relies on a parameter-free and
problem-independent Artificial Neural Network based (ANN) approach to
discontinuity detection and assignment of localized viscosity
$\mu[\textbf{e}]$, that was introduced
in~\cite{schwander2021controlling} and was subsequently modified and
extended to the non-periodic FC-based context
in~\cite{bruno2022fc}. Briefly, the approach relies on a single
pre-trained neural network for detection of discontinuities on the
basis of Gibbs oscillations in Fourier series, allowing for a sharp
resolution of flow features such as shocks and contact
discontinuities.  A description of the artificial viscosity-assignment
method utilized in this paper, including techniques that, as suggested
above, are necessary in the present general-domain multi-patch
context, is presented in Section~\ref{sec:viscosity}.  On the basis of
such discontinuity-smearing approaches, the FC-SDNN method produces
smooth flows away from discontinuities, sharp shocks and contact
discontinuities, and accurate flow predictions that converge rapidly
as the underlying discretizations are suitably refined.

\section{\label{sec:Multidomain} Overlapping-patch/subpatch geometry description and discretization}
\subsection{\label{subsec:decomposition} Domain decomposition: patches and their parametrizations}

\begin{figure}[H]
\centering
\includegraphics[width=0.6\linewidth,]{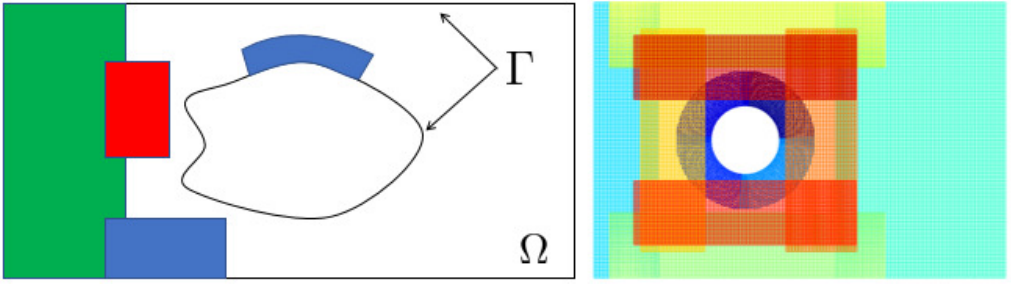}
\caption{Left: Patch types used in the overlapping-patch decomposition
  of a general domain $\Omega$, including, (a)~A corner patch of type
  $\mathcal{C}$ (green); (b)~Smooth-boundary patches of type
  $\mathcal{S}$ (blue); (c)~An interior patch of type $\mathcal{I}$
  (red). Right: full patch decomposition of a computational domain
  containing a circular obstacle.}
    \label{fig:mappings}
\end{figure}

We seek to produce an FC-based discretization of equation~\eqref{eq:
  euler 2d equation} over a general 2D domain $\Omega$ delimited by
obstacles with smooth boundaries together with an outer computational
domain boundary, as depicted in Figure~\ref{fig:mappings}. We denote
by $\Gamma = \partial \Omega$ the boundary of the complete
computational domain, which includes the smooth obstacle boundary or
boundaries as well as the outer computational domain boundary---the
latter one of which is throughout this paper assumed to be rectangular
in shape.
Following~\cite{chesshire1990composite,albin2011spectral,amlani2016fc,bruno2019higher},
in order to obtain the desired FC discretization the proposed method
relies on the decomposition of $\Omega$ as a union of a number
$P = P_{\mathcal{S}} + P_{\mathcal{C}}+ P_{\mathcal{I}}$ of {\em
  overlapping} patches, each one of which is an open set, including, a
number $P_{\mathcal{S}}$ of $\mathcal{S}$-type patches
$\Omega^{\mathcal{S}}_{p}$ (smooth-boundary patches); a number
$P_{\mathcal{C}}$ of $\mathcal{C}$-type patches (corner patches)
$ \Omega^{\mathcal{C}}_{p}$; and a number $P_{\mathcal{I}}$ of
$\mathcal{I}$-type patches $\Omega^{\mathcal{I}}_{p}$ (interior
patches):
\begin{equation}\label{decomp}
  \Omega = \left(\bigcup_{p=1}^{P_{\mathcal{S}}} \Omega^{\mathcal{S}}_{p}\right ) \cup \left(\bigcup_{p=1}^{P_{\mathcal{C}}} \Omega^{\mathcal{C}}_{p}\right )  \cup \left(\bigcup_{p=1}^{P_{\mathcal{I}}} \Omega^{\mathcal{I}}_{p}\right ).
\end{equation}
see also Remark~\ref{rem_bdries}.  Here $\mathcal{C}$-type patches
(shown in green in Figure~\ref{fig:mappings}) contain a neighborhood
within $\Omega$ of one or more boundary corner points in the
rectangular outer boundary; $\mathcal{S}$-type patches (shown in blue
in the figure) cover areas along smooth portions of the outer boundary
and portions of the boundary of the obstacle(s) (which are assumed to
be smooth throughout this paper, but however
see~\cite{bruno_leibo_partII} where non-smooth obstacle boundaries are
allowed); and, $\mathcal{I}$-type patches (shown in red in the figure)
cover regions away from boundaries. As detailed in what follows, each
patch $\Omega^{\mathcal{R}}_{p}$ ($\mathcal{R}=\mathcal{S}$,
$\mathcal{I}$, $\mathcal{C}$, $p=1,\dots,P_{\mathcal{R}}$) is
constructed and discretized on the basis of a smooth invertible
parametrization
\begin{equation}\label{eqn:mappings}
 \mathcal{M}^\mathcal{R}_{p}: \mathcal{Q} \rightarrow \overline{\Omega^{\mathcal{R}}_{p}}, \qquad
    \mathcal{Q} \coloneqq [0,1]\times[0,1], 
\end{equation}
which maps the canonical (square) parameter space $\mathcal{Q}$ onto
the {\em closure} $\overline{\Omega^{\mathcal{R}}_{p}}$ of the patch.
\begin{remark}\label{rem_bdries}
  A requirement imposed on the patches $\Omega^{\mathcal{R}}_{p}$,
  which is satisfied in all cases by the patches constructed in the
  following sections, is that the image of each of the sides of the
  parameter domain $\mathcal{Q}$ is either fully contained within
  $\Gamma$, or it intersects $\Gamma$ at exactly one point (which
  occurs as the image of the relevant side of $\mathcal{Q}$ meets
  $\Gamma$ transversally), or its intersection with $\Gamma$ is empty.
\end{remark}
Using the patch parametrizations~\eqref{eqn:mappings} the necessary
discretization of each patch $\Omega^{\mathcal{R}}_{p}$ is produced as
the image by the mapping $\mathcal{M}^\mathcal{R}_{p}$ of a simple
Cartesian discretization of $\mathcal{Q}$---with possibly different
meshsizes along the $q^1$ (abscissa) and $q^2$ (ordinate) directions
in $\mathcal{Q}$.  The required overlaps between patches are assumed
to be ``sufficiently broad''---so that, roughly speaking, except for
discretization points near physical boundaries, each discretization
point $\mathbf{x}$ associated with a given patch
$\Omega^{\mathcal{R}}_p$, ($\mathcal{R}=\mathcal{S}$, $\mathcal{I}$,
$\mathcal{C}$; $p=1,\dots,P_{\mathcal{R}}$), and located in the
vicinity of the boundary of this patch, is contained in a
``sufficiently interior'' region of at least one other patch, say
$\Omega^{\mathcal{R}}_{p'}$
($(\mathcal{R}', p')\ne (\mathcal{R}, p)$). This overlap condition is
imposed so as to ensure adequate interpolation ranges between patches
as well as creation of sufficiently smooth viscosity windows near
patch boundaries---as detailed in Sections~\ref{subsec:overlap}
and~\ref{subsec:window_propagation}, respectively.  To enable the use
of arbitrarily fine meshes while avoiding Fourier expansions of
excessively high order—thereby maintaining control over
differentiation errors, especially near sharp flow features such as
shock waves and other solution discontinuities—and to further
facilitate parallelization, the proposed algorithm permits
partitioning of the patches $\Omega^{\mathcal{R}}_p$
($\mathcal{R}=\mathcal{S}, \mathcal{I}, \mathcal{C}$) into sets of
similarly overlapping ``subpatches,'' which are introduced in
Section~\ref{subsec:overl_dec}.

Clearly, $\mathcal{C}$- and $\mathcal{I}$-type patches can be
constructed on the basis of affine parametrizations
$\mathcal{M}^{\mathcal{R}}_p:\mathcal{Q}\to\overline{\Omega^{\mathcal{R}}_p}$
($\mathcal{R} = \mathcal{I}$ or $\mathcal{R} = \mathcal{C}$). The
$\mathcal{S}$-type patches adjacent to smooth segments of $\Gamma$, in
turn, can be constructed easily. Indeed, consider a smooth arc
$\overset{\frown}{AB} \subset \Gamma$ connecting two selected points
$A, B \in \Gamma$, and let $\ell_A:[0,1] \to \mathbb{R}^2$ be a
smooth, invertible parameterization of this arc such that
\begin{equation}\label{corner_mapping}
  \ell_A([0, 1]) = \overset{\frown}{AB}.
\end{equation}
Letting $\nu:[0,1]\to\mathbb{R}^2$ denote the unit normal vector to
$\overset{\frown}{AB}$ (so that, for $t\in[0,1]$, $\nu(t)$ is the unit
normal vector at the point $\ell_A(t)$ pointing towards the interior
of $\Omega$), and letting $H>0$ denote an adequately selected
parameter, the desired 2D $\mathcal{S}$-type patch parametrization
$\mathcal{M}^{\mathcal{S}}_p:\mathcal{Q}\to\overline{\Omega^{\mathcal{S}}_p}$,
illustrated in Figure~\ref{fig:mappings}, can be produced via the
expression
\begin{equation}
  \label{eq:param_S}
  \mathcal{M}^{\mathcal{S}}_p(q^1,q^2) = \ell_A(q^1)+ q^2 H  \nu(q^1).
\end{equation}

\subsection{\label{subsec:overl_dec}Patches, subpatches and their grids}

An overall overlapping-patch decomposition of the form~\eqref{decomp}
for the domain $\Omega$ is obtained on the basis of adequately
selected sets of patches, including patches of types $\mathcal{C}$,
$\mathcal{S}$ and $\mathcal{I}$, as described in
Section~\ref{subsec:decomposition}. Since $\Omega$ is a connected
open set, and since all patches are themselves open sets, each patch
$\Omega^{\mathcal{R}}_p$ ($\mathcal{R}=\mathcal{S}$, $\mathcal{C}$, or
$\mathcal{I}$ and $1\leq p\leq P_\mathcal{R}$) must necessarily
overlap one or more patches $\Omega^{\mathcal{R}'}_{p'}$
($\mathcal{R}'=\mathcal{S}$, $\mathcal{C}$, or $\mathcal{I}$ and
$1\leq p'\leq P_{\mathcal{R}'}$). In fact, an efficient implementation
of the proposed algorithm additionally requires that a ``sufficient
amount of overlap'' exists, in the sense that every point
$x\in\Omega^{\mathcal{R}}_p$ that lies ``in the vicinity'' of the
boundary of $\Omega^{\mathcal{R}}_p$ must lie ``sufficiently deep''
within some patch $\Omega^{\mathcal{R}'}_{p'}$ (i.e., within
$\Omega^{\mathcal{R}'}_{p'}$ and away from a vicinity of the boundary
of $\Omega^{\mathcal{R}'}_{p'}$) for some
$(p',\mathcal{R}')\ne (p,\mathcal{R})$.  The vicinity and depth
concepts alluded to above are defined and quantified in
Sections~\ref{subsec:overlap}, \ref{subsubsec:interpatch}
and~\ref{subsubsec:intrapatch} on the basis of the patch
discretizations used---which are themselves introduced in
Section~\ref{square}.

For flexibility in both geometrical representation and discretization
refinement, the patching structure we use incorporates a decomposition
in ``sub-patches'' $\Omega^{\mathcal{R}}_{p, \ell}$ of each patch
$\Omega^{\mathcal{R}}_{p}$ ($\mathcal{R}=\mathcal{S}$,
$\mathcal{C}$, or $\mathcal{I}$):
\[
  \Omega^{\mathcal{R}}_{p} = \bigcup_{\ell=1}^{r_p} \Omega^{\mathcal{R}}_{p, \ell},
\]
on each of which a patch-dependent number of discretization points is
used, as indicated in Section~\ref{square}.  Like the patches
themselves, the subpatches we use are required to display sufficiently
large amounts of overlap with their neighbors. Subpatches are defined,
simply, via overlapping decomposition of the canonical parameter space
$\mathcal{Q}$, as described in Section~\ref{square}. The description
in that section utilizes detailed notations for the discretizations
used for the various patch and subpatch types. Such level of detail,
which involves several sub- and super-indexes, is needed in the
context of the description of the multi-patch viscosity assignment and
propagation strategy presented in Section~\ref{sec:viscosity}.

\subsubsection{ $\mathcal{Q}$ parameter space discretization for $\mathcal{S}$-, $\mathcal{I}$- and $\mathcal{C}$-type
  patches
  \label{square}}

\begin{figure}[H]
\centering
    \includegraphics[width=0.6\linewidth,]{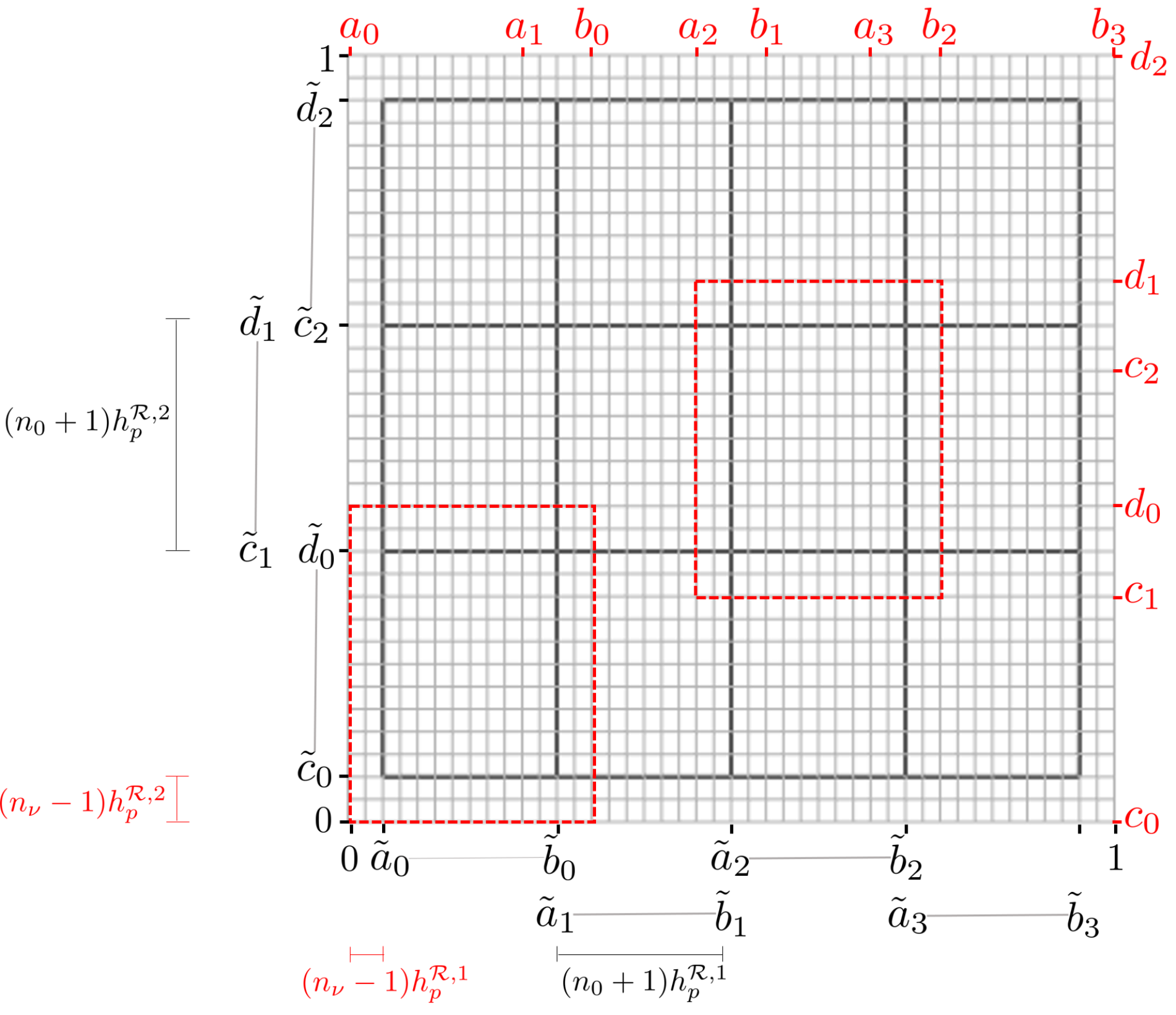}
    \caption{Patch decomposition (in $\mathcal{Q}$ parameter-space)
      for a patch $\Omega^{\mathcal{R}}_p$ (with
      $\mathcal{R}=\mathcal{S}$, $\mathcal{C}$, or $\mathcal{I}$)
      into preliminary non-overlapping rectangular subpatches (with
      boundaries shown in black), associated overlapping rectangles
      (two examples highlighted, with boundaries shown in red) and
      discretization lines (shown in gray). Geometric parameter values
      $r_p = 4$, $s_p = 3$, $n_0 = 9$ and $n_v = 3$ were used for
      this illustration.}
    \label{fig:qpatch}
\end{figure}

In order to obtain discretizations for each one of the patches
$\Omega^\mathcal{R}_p$ ($\mathcal{R}=\mathcal{S}$, $\mathcal{C}$,
$\mathcal{I}$; $1\leq p\leq P_\mathcal{R}$), we introduce the
parameter space grids
\begin{equation}
  \label{eq:param-grid}
  G^{\mathcal{R}}_{p} = \mathcal{Q} \cap \big\{(q^{\mathcal{R}, 1}_{p, i}, q^{\mathcal{R}, 2}_{p, j})\ :\ q^{\mathcal{R}, 1}_{p, i} = ih^{\mathcal{R}, 1}_{p}, q^{\mathcal{R}, 2}_{p, j} =
  jh^{\mathcal{R}, 2}_p, \quad (i, j) \in \mathbb{N}_0^2 = \mathbb{N}\cup\{ 0\} \big\},
\end{equation}
of grid sizes
$h^{\mathcal{R}, 1}_{p} = \frac{1}{N^{\mathcal{R}, 1}_p + 2n_v - 1}$
and
$h^{\mathcal{R}, 2}_{p} = \frac{1}{N^{\mathcal{R}, 2}_p + 2n_v - 1}$
along the $q^1$ (abscissa) and $q^2$ (ordinate) directions in
parameter space, respectively. Here $N^{\mathcal{R}, 1}_p$ (resp.
$N^{\mathcal{R}, 2}_p$) denotes the number of discretization points,
in each discretization line along the $q^1$ (resp. $q^2$) direction
within the square $\mathcal{Q}$ (equation~\eqref{eqn:mappings}), that lie outside
a certain layer of discretization points in a ``boundary vicinity''
region near the boundary of $\mathcal{Q}$. The boundary vicinity is
itself characterized by the positive integer $n_v$ equal to the number
of discretization points along the boundary-vicinity width, as
illustrated in Figure~\ref{fig:qpatch}.  The value $n_v = 9$ was used
in all of the implementations presented in this paper.

As indicated in Section~\ref{subsec:decomposition}, the proposed
discretization algorithm additionally utilizes subpatches to limit the
order of the FC expansions used. To describe the sub-patch
decomposition approach we first introduce, for a given assignment of
positive integers $r^{\mathcal{R}}_p$ and $s^{\mathcal{R}}_p$ for each
$(\mathcal{R},p)$, a set of non-overlapping rectangles strictly
contained within $\mathcal{Q}$:
\begin{equation}
  \label{eq:rect}
  \bigcup_{(r, s) \in \Theta^{\mathcal{R}}_{p}} [\tilde{a}_{r}, \tilde{b}_{r}] \times [\tilde{c}_{s}, \tilde{d}_{s}] \subsetneq \mathcal{Q}
\end{equation}
where we define
\begin{equation}
  \label{eq:sp_index}
  \Theta^{\mathcal{R}}_{p} = \{(r,s)\in\mathbb{N}_0^2\ |\ 0\leq r\leq r^{\mathcal{R}}_p
  - 1 \mbox{ and } 0\leq s\leq s^{\mathcal{R}}_p - 1 \}.
\end{equation}
In what follows we select
$N^{\mathcal{R}, 1}_p = r^{\mathcal{R}}_p(n_0 + 1) - 1$ and
$N^{\mathcal{R}, 2}_p = s^{\mathcal{R}}_p(n_0 + 1) - 1$, where
$r^{\mathcal{R}}_p$ and $s^{\mathcal{R}}_p$ are given positive
integers, and where $n_0$ denotes the number of points used in the
discretization of certain preliminary non-overlapping subpatches
introduced below along both the parameter space $q^1$ and $q^2$
directions. The procedure used for selection of the integers
$r^{\mathcal{R}}_p$ and $s^{\mathcal{R}}_p$ is described in
Section~\ref{subsubsec:subpatches}. The $(\mathcal{R}, p)$-dependent
interval endpoints $\tilde{a}_{r}$, $\tilde{b}_{r}$, $\tilde{c}_{s}$
and $\tilde{d}_{s}$ are selected in terms of the integer parameters
$n_v$ and $n_0$ introduced above.

\begin{equation}\label{discretization}
  \begin{aligned}
   & \tilde{a}_0 = (n_{v}-1) h^{\mathcal{R}, 1}_p, &\\
   & \tilde{b}_{r} - \tilde{a}_r = (n_0 + 1) h^{\mathcal{R}, 1}_p, &\qquad (0 \leq r \leq r^{\mathcal{R}}_p - 1),  \\
   & \tilde{a}_{r+1} = \tilde{b}_{r},  &\qquad (0 \leq r < r^{\mathcal{R}}_p - 1),  \\
   & \tilde{c}_0 = (n_{v}-1) h^{\mathcal{R}, 2}_p, & \\
   & \tilde{d}_{s} - \tilde{c}_s = (n_0+1)  h^{\mathcal{R}, 2}_p,  &\qquad (0 \leq s \leq s^{\mathcal{R}}_p - 1 ), \\
   & \tilde{c}_{s+1} = \tilde{d}_{s},  &\qquad (0 \leq s < s^{\mathcal{R}}_p - 1).
  \end{aligned}
\end{equation}
  In all of the implementations presented in this paper the value
  $n_0 = 83$ was used.

   The overlapping subpatches associated with a given patch
  are characterized, via the patch parametrization, by the
  decomposition of $\mathcal{Q}$ as the union
\begin{equation}
  \label{eq:Spart}
  \mathcal{Q} =  \bigcup_{(r, s) \in \Theta^{\mathcal{R}}_{p}} [a_{r}, b_{r}] \times [c_{s}, d_{s}],
\end{equation}
of overlapping rectangles (Figure~\ref{fig:qpatch}), whose vertices
are given by

\begin{equation}\label{subpatchesQ}
  \begin{aligned}
    &a_{r}  = \tilde{a}_r - (n_{v}-1) h^{\mathcal{R}, 1}_p, &\qquad (0 \leq r \leq r^{\mathcal{R}}_p - 1),  \\
    &b_{r}  = \tilde{b}_{r} + (n_{v}-1) h^{\mathcal{R}, 1}_p, &\qquad (0 \leq r \leq r^{\mathcal{R}}_p - 1), \\
    &c_{s}  = \tilde{c}_s - (n_{v}-1) h^{\mathcal{R}, 2}_p, &\qquad (0 \leq s \leq s^{\mathcal{R}}_p - 1),  \\
    &d_{s}  = \tilde{d}_{s} + (n_{v}-1) h^{\mathcal{R}, 2}_p, &\qquad (0 \leq s \leq s^{\mathcal{R}}_p - 1).  
  \end{aligned}
\end{equation}
Note that each subpatch contains a total of
\begin{equation}
  \label{eq:pts_p_subp}
  N_\mathcal{Q}= (n_0+2n_v)^2
\end{equation}
discretization points.

For $\ell = (r, s) \in \Theta^{\mathcal{R}}_p$, we denote by
$G^{\mathcal{R}}_{p, \ell}$ the parameter space grid for the $\ell$-th
subpatch of the patch $\Omega^{\mathcal{R}}_p$; we clearly have 
\[ G^{\mathcal{R}}_{p, \ell} = G^{\mathcal{R}}_{p} \cap \left( [a_{r},
    b_{r}] \times [c_{s}, d_{s}] \right).
\]
Using the mappings $\mathcal{M}^{\mathcal{R}}_p$
($\mathcal{R}=\mathcal{C}$, $\mathcal{S}$, $\mathcal{I}$) described in
Section~\ref{subsec:decomposition}, we obtain the desired subpatches
$\Omega^{\mathcal{R}}_{p, \ell} = \mathcal{M}^{\mathcal{R}}_p([a_r,
b_r]\times[c_s, d_s])$ ($\ell = (r,s)$) and the physical grids
\begin{equation}
  \label{eq:3_grids}
  \mathcal{G}^{\mathcal{R}}_p = \mathcal{M}^{\mathcal{R}}_p\left( G^{\mathcal{R}}_p \right) \quad \mbox{and} \quad \mathcal{G}^{\mathcal{R}}_{p, \ell} = \mathcal{M}^{\mathcal{R}}_p\left( G^{\mathcal{R}}_{p, \ell}\right).
\end{equation}
The set of indices $\mathcal{D}^{\mathcal{R}}_{p, \ell}$ of parameter
space grid points contained in $G^{\mathcal{R}}_{p, \ell}$, finally,
is denoted by
\begin{equation} \label{D_R}
  \mathcal{D}^{\mathcal{R}}_{p, \ell} = \{(i, j) \in \mathbb{N}_0^2 \ |\ (ih^{\mathcal{R}, 1}_p, jh^{\mathcal{R}, 2}_p) \in G^{\mathcal{R}}_{p, \ell} \}, \quad
  \mathcal{R} = \mathcal{C}, \mathcal{S}, \mathcal{I}.
\end{equation}

\subsubsection{\label{subsec:overlap}Minimum subpatch overlap
  condition}

As indicated in Sections~\ref{subsubsec:interpatch}
and~\ref{subsubsec:intrapatch} below, and as alluded to previously in
Section~\ref{subsec:decomposition}, in order to properly enable
inter-patch data communication and the overall multi-patch viscosity
assignment, the overlap between patches $\Omega^{\mathcal{R}}_p$ must
be sufficiently broad.  The required overlap breadth is quantified in
terms of a parameter associated with the selected subpatches of each
patch (namely, the integer parameter $n_v$ introduced in
Section~\ref{square}) as well as the patch-boundary ``sides''. (The
``sides'' of a patch $\Omega^{\mathcal{R}}_p$ or a subpatch
$\Omega^{\mathcal{R}}_{p, \ell}$ are defined as the images under the
corresponding patch parametrization of the straight segments that make
up the boundary of the corresponding parameter rectangle, that is, the
sides of the square $\mathcal{Q}$ for $\Omega^{\mathcal{R}}_p$, and
the sides of the relevant rectangle in equation~\eqref{eq:Spart} for
$\Omega^{\mathcal{R}}_{p, \ell}$.)  Utilizing these concepts, the
minimum-overlap condition is said to be satisfied for a given patch
$\Omega^{\mathcal{R}}_p$ if and only if for each side of
$\Omega^{\mathcal{R}}_p$ that is not contained in $\partial \Omega$
(cf. Remark~\ref{rem_bdries}), a $(2 n_v + 1)$-point wide layer of
grid points adjacent to that side in the $\Omega^{\mathcal{R}}_p$ grid
is also included in a union of one or more patches
$\Omega^{\mathcal{R'}}_{p'}$ with
$(p',\mathcal{R}')\neq (p,\mathcal{R})$. Note that, per the
construction in Section~\ref{square} at the parameter space level,
subpatches of a single patch satisfy the minimum overlap condition
embodied in equation~\eqref{subpatchesQ}: neighboring subpatches share
a $(2 n_v + 1)$-point wide overlap.

\subsubsection{Subpatches meshsizes and grid refinement\label{subsubsec:subpatches}}

Very general overlapping patch decompositions, satisfying the minimum
overlap condition introduced in Section~\ref{subsec:overlap}, can be
obtained on the basis of the procedures described in
Section~\ref{subsec:decomposition}. As briefly discussed below in this
section, further, the strategies presented in Section~\ref{square} can
be used to produce sets of overlapping subpatches, each one endowed
with a subpatch grid, in such a way that, 1)~The set of all patches
(namely, all $\mathcal{S}$-, $\mathcal{C}$- and $\mathcal{I}$-type
patches) satisfies the overlap conditions introduced in
Section~\ref{subsec:overlap}; 2)~Each subpatch contains exactly
$N_\mathcal{Q}=\left(n_0 + 2 n_v\right)^2$ discretization points; and,
3)~All of the subpatch parameter-space grid sizes
$h^{\mathcal{R}, 1}_p$ and $h^{\mathcal{R}, 2}_p$
($\mathcal{R}=\mathcal{S}$, $\mathcal{C}$, $\mathcal{I}$;
$1\leq p\leq P_\mathcal{R}$) are such that the resulting physical grid
sizes (defined as the maximum distance between two consecutive grid
points in physical space) are less than or equal to a user-prescribed
upper bound $\overline{h}>0$.

Indeed, as indicated in Sections~\ref{square}
and~\ref{subsec:overlap}, the sets of patches and subpatches
constructed per the description in Section~\ref{square} satisfy
points~1) and~2), but, clearly, they may or may not satisfy point~3)
for certain values of $\mathcal{R}$ and $p$. For the values of
$\mathcal{R}$ and $p$ that do not satisfy point~3) the integers
$r^{\mathcal{R}}_p$ and $ s^{\mathcal{R}}_p$ (cf.
Section~\ref{square}) are suitably increased---thus proportionally
increasing the numbers of subpatches and the number of gridpoints of
the patch $\Omega^{\mathcal{R}}_p$ while leaving the overall patch
decomposition and mappings $\mathcal{M}^{\mathcal{R}}_p$
unchanged---until the physical grid size upper bound condition in
point~3) is satisfied, as desired.

Once a set of patches, subpatches and associated discretizations has
been obtained in such a way that points~1) through~3) are satisfied
for a given $\overline{h}$ value, say $\overline{h} = \overline{h}_0$,
additional global grid refinements may be obtained by multiplying the
integers $r^{\mathcal{R}}_p$ and $ s^{\mathcal{R}}_p$ for all
$\mathcal{R}$ and $p$ values by a constant integer factor $K$---which
results in a mesh satisfying condition 3) with $\overline{h}$
approximately equal to $\frac{\overline{h}_0}{K}$. Naturally, the
overall grid can also be fine-tuned by locally adjusting the mesh size
on one or more patches $\Omega^{\mathcal{R}}_p$, in either or both
coordinate directions $q^1$ and $q^2$, through appropriate
modification of the integers $r^{\mathcal{R}}_p$ and
$s^{\mathcal{R}}_p$ associated with each patch.

\subsection{\label{subsec:curvilinear_eqns}Governing equations in curvilinear coordinates}

In the multi-patch setting described in this section, the solution
vector $\textbf{e}(\textbf{x}, t)$ can be viewed as a family
of vectors
$\textbf{e}^{\mathcal{R}}_{p, \ell}(\textbf{x}, t)$ (with
$\mathcal{R}=\mathcal{S}$, $\mathcal{C}$, or
$\mathcal{I}$, $1\leq p\leq P_\mathcal{R}$ and
$\ell \in \Theta^{\mathcal{R}}_p$), defined on
$\Omega^{\mathcal{R}}_{p, \ell} \times \mathbb{R}^{+}$.  To evolve the
solution on the curvilinear grids
$\mathcal{G}^{\mathcal{R}}_{p, \ell}$ associated to the smooth
invertible mapping $\mathcal{M}^{\mathcal{R}}_{p}$, the Jacobian
matrix $J^{\mathcal{R}, p}_{\mathbf{q}}(\mathbf{x})$ of the inverse
mapping $\big(\mathcal{M}^{\mathcal{R}}_{p}\big)^{-1}$ is
utilized. Thus, using the chain rule expression
$\nabla_{\mathbf{x}} = [ J^{\mathcal{R},
  p}_{\mathbf{q}}(\mathbf{x})]^T \nabla_{\mathbf{q}}$, we obtain the curvilinear
form
\begin{equation} \label{eq:general_curvilinear}
  \frac{\partial \textbf{e}^{\mathcal{R}}_{p, \ell}(\textbf{x}, t)}{\partial t} + \big[J^{\mathcal{R}, p}_{\textbf{q}} (\textbf{x})\big]^T \nabla_{\textbf{q}} \cdot
  \big( f(\textbf{e}^{\mathcal{R}}_{p, \ell}(\textbf{x}, t))\big) = \big[J^{\mathcal{R}, p}_{\textbf{q}} (\textbf{x})\big]^T \nabla_{\textbf{q}} \cdot
  \left( f_{visc}\big(\textbf{e}^{\mathcal{R}}_{p, \ell}(\textbf{x}, t)\big)\right), \quad \textbf{x} \in \Omega^{\mathcal{R}}_{p, \ell}.
\end{equation}
of the artificial viscosity-augmented hyperbolic system~\eqref{eq: convection diffusion eqn}.

\subsection{\label{subsec:connectivity} Patch/subpatch communication
  of solution values}

\begin{figure}[H]
\centering
    \includegraphics[width=1\linewidth,]{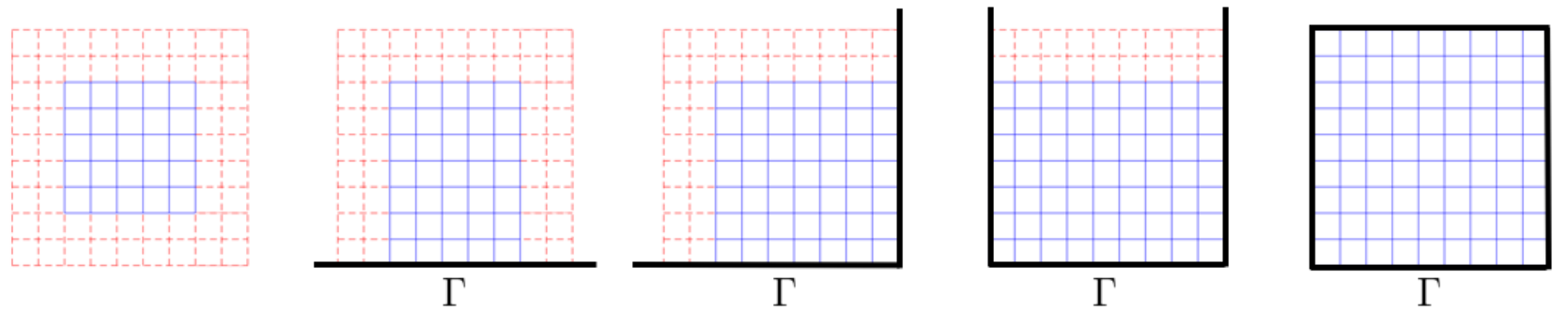}
    \caption{Illustration of subpatches (union of blue and red) and
      their fringes (red) for various kinds of subpatches. For
      simplicity this illustration uses the value $n_f = 2$ ($n_f = 5$
      was used for all the numerical examples presented in this
      paper). Left panel: all the sides of the subpatch are internal
      sides. Four rightmost panels, from left to right: subpatches
      containing one, two, three and all four external sides, all of
      which, by design, are contained in the domain boundary $\Gamma$
      (black).}
    \label{fig:fringes}
\end{figure}

The overlap built into the patch/subpatch domain decomposition
approach described in the previous sections enables the communication
(possibly via interpolation) of subpatch grid-point values of the
solution $\textbf{e}$ in the vicinity of boundaries of
subpatches---which is a necessary element in the overlapping-patch
time-stepping algorithm introduced in
Section~\ref{subsec:time_marching}.  The communication of solution
values relies on the concept of ``$n$-point fringe region'' of a
subpatch.  To introduce this concept we consider the sides of any
subpatch $\Omega^{\mathcal{R}}_{p, \ell}$ (defined in
Section~\ref{subsec:overlap}), and we note that, per the constructions
in that section, each side of a subpatch is either contained within
$\Gamma$ or it intersects $\Gamma$ at most at one point. In the first
(resp. the second) case we say that the side is ``external''
(resp. ``internal''). Then, denoting by
$\mathcal{I}^{\mathcal{R}}_{ p, \ell}$ the set of all grid points in
$\mathcal{G}^{\mathcal{R}}_{p, \ell}$ that are contained in the
internal sides of $\Omega^{\mathcal{R}}_{p, \ell}$, we define the
$n$-point fringe region $\mathcal{F}^{\mathcal{R}}_{p, \ell, n}$
of $\Omega^{\mathcal{R}}_{p, \ell}$ as the set of all grid points
$\mathbf{x}_{ij}= (x_i,y_j)\in \mathcal{G}^{\mathcal{R}}_{p, \ell}$
whose distance to $\mathcal{I}^{\mathcal{R}}_{ p, \ell}$ (in the sense
of the ``maximum index norm'' distance $d$ defined below) is less than
or equal to $n$. In detail, the fringe region is given by
\begin{equation}\label{fringe_region}
 \mathcal{F}^{\mathcal{R}}_{p, \ell, n} \coloneqq \{ \mathbf{x}_{ij} \in \mathcal{G}^{\mathcal{R}}_{p, \ell}, \quad d(\mathbf{x}_{ij}, \mathcal{I}^{\mathcal{R}}_{ p, \ell}) < n\}
\end{equation}
where the maximum index norm distance
$ d(\mathbf{x}_{ij}, \mathcal{I}^{\mathcal{R}}_{ p, \ell})$ from the
point $\mathbf{x}_{ij} \in \mathcal{G}^{\mathcal{R}}_{p, \ell}$ to the set
$\mathcal{I}^{\mathcal{R}}_{p, \ell}$ is defined by
\begin{equation}
   d(\mathbf{x}_{ij}, \mathcal{I}^{\mathcal{R}}_{ p, \ell}) \coloneqq  \min_{\mathbf{x}_{rs} \in \mathcal{I}^{\mathcal{R}}_{p, \ell}} \max \left\{ |i - r|, |j - s| \right\}.
\end{equation}
Figure~\ref{fig:fringes} illustrates the different fringe regions that
can arise within a subpatch, depending on the internal sides contained
in the subpatch. Additionally, fringe regions with two different
values of $n$, namely, $n=n_f$ and $n=n_v$, are used by the algorithm to
tackle the problems of inter- and intra-patch data communication and
viscosity assignment, respectively, as described in what follows. The
values $n_f=5$ and $n_v=9$ were used for all of the numerical examples
presented in both parts, Part~I and Part~II, of this two-part
contribution.

As described over the next few sections, as part of the time-stepping
process each subpatch $\Omega^{\mathcal{R}}_{p, \ell}$ receives
solution data (from neighboring subpatches) on its $n$-fringe region
$\mathcal{F}^{\mathcal{R}}_{p, \ell, n_f}$.  In fact, the
time-stepping method we use utilizes two different algorithms for the
communication of solution data between pairs of subpatches, including
``inter-patch'' communication (that is, communication between two
subpatches of different patches, which requires polynomial {\em
  interpolation} of grid-point values of the solution $\textbf{e}$)
and ``intra-patch'' communication (which only involves {\em exchange}
of grid-point values of $\textbf{e}$ between two subpatches of the
same patch). Details concerning these two procedures are presented in
Sections~\ref{subsubsec:interpatch} and~\ref{subsubsec:intrapatch},
respectively. As detailed in the following remark, the minimum overlap
condition is required to ensure that donor points may be selected
which are not recipients of solution data.

\begin{remark}\label{overlap_cond}
  The overlapping patch decomposition and discretization should be set
  up in such a way that data donor grid points are not themselves
  receivers of data from other patches. If this condition is not
  satisfied then the well-posedness of the problem and, thus,
  stability of the time stepping method, are
  compromised~\cite{chesshire1990composite}.  This donor-receiver
  condition is generally satisfied by requiring that donor grid points
  are located sufficiently deeply within the subpatch, while receiving
  grid points are located next to the boundary of the subpatch, as
  detailed in Section~\ref{subsubsec:interpatch} for the case of
  solution-value communication between subpatches of different
  patches, and in Section~\ref{subsubsec:intrapatch} for the case of
  communication between subpatches of a single patch. Additionally,
  use of donor data from regions well separated from subpatch
  boundaries is advantageous since the FC representation of the
  solution (an indeed, the numerical representation resulting from any
  discretization method) is more accurate and smoother away from the
  boundary of the representation region---which in this case are the
  subpatch boundaries themselves.
\end{remark}

\subsubsection{\label{subsubsec:interpatch} Inter-patch data communication}

\begin{figure}[H]
\centering
    \includegraphics[width=0.5\linewidth,]{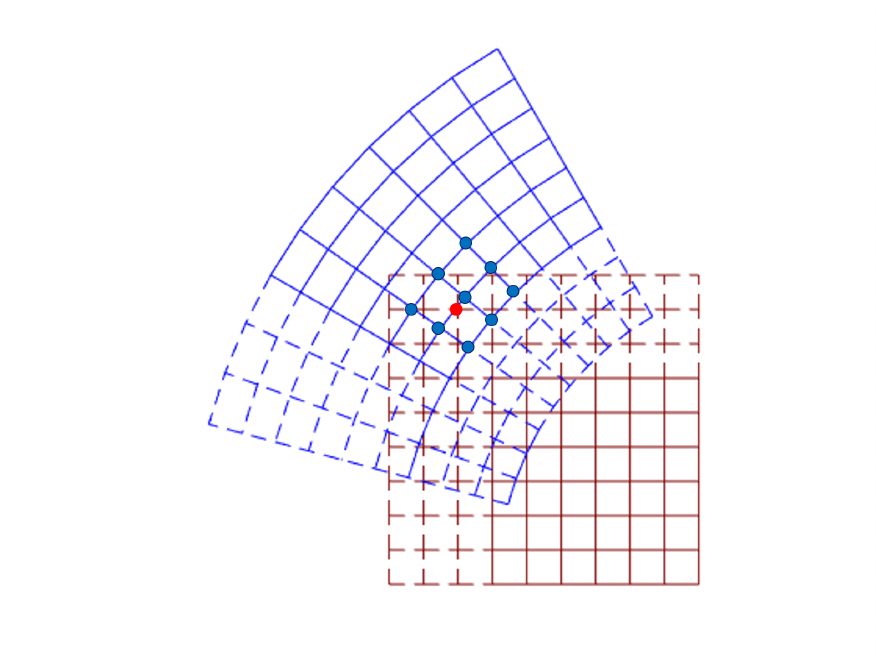}
    \caption{Illustration (with $n_f = 3$, for graphical simplicity)
      of the inter-patch data communication effected via exchanges
      between the blue and brown subpatches (which, as befits
      inter-patch communications, are subpatches of different
      patches). Dashed-line regions: receiving fringe regions for each
      subpatch. Red dot: a receiver point in the fringe region of the
      brown subpatch, set to receive solution data from the blue
      subpatch, via interpolation. Blue dots: the $3\times3$
      donor-point stencil in the blue subpatch used for interpolation
      to the red receiver point. As indicated in
      Remark~\ref{overlap_cond}, none of the points in the donor
      stencil are located in the fringe region of the blue subpatch.}
    \label{fig:interpatch}
  \end{figure}

As indicated above, the proposed time-stepping algorithm requires
communication of values of the solution vector
$\textbf{e}= \textbf{e}(\mathbf{x},t)$ at grid points in certain
boundary-vicinity regions of each subpatch
$\Omega^{\mathcal{R}}_{p,\ell}$; in this section we consider only
communication between subpatches $\Omega^{\mathcal{R}}_{p,\ell}$ and
$\Omega^{\mathcal{R}'}_{p',\ell'}$ with
$(\mathcal{R}', p') \neq (\mathcal{R}, p)$ (that is, between
subpatches of different patches). Thus, prior to every time step (or,
more precisely, prior to each time-stage in the SSPRK-4 time-stepping
method we use, see Section~\ref{subsec:differentiation}), for each
$(\mathcal{R},p,\ell)$ and for every point $\mathbf{x}$ in the fringe
region $\mathcal{F}^{\mathcal{R}}_{p, \ell, n_f}$ of the subpatch
$\Omega^{\mathcal{R}}_{p,\ell}$, the solution value
$\textbf{e}(\mathbf{x})$ is overwritten with the value obtained for
the same quantity by interpolation from an adjacent donor subpatch
$\Omega^{\mathcal{R}'}_{p',\ell'}$. (Among the possibly multiple
neighboring donor subpatches $\Omega^{\mathcal{R}'}_{p',\ell'}$
containing $\mathbf{x}$, with
$(\mathcal{R}', p') \neq (\mathcal{R}, p)$, a donor subpatch is
selected for which $\mathbf{x}$ is the farthest from the boundary of
$(\mathcal{R}', p')$.)  A $3 \times 3$-point stencil of
$\mathcal{G}^{\mathcal{R}'}_{p',\ell'}$ grid points surrounding
$\mathbf{x}$, together with an interpolation scheme based on iterated
1D second order polynomial
interpolation~\cite[Sec. 3.6]{press2007numerical} in parameter space,
enacted on the basis of the Neville algorithm, is used to produce the
necessary interpolated values, as illustrated in
Figure~\ref{fig:interpatch}.

With reference to Remark~\ref{overlap_cond}, the algorithm developed
in this paper includes a testing phase which, for a given
patch/subpatch decomposition and grid assignments, checks whether, as
indicated in that remark, any data donor points in any patch are
themselves receivers of solution data. If this donor-receiver
condition is not satisfied, possible remedies include 1)~Refining the
grids by proportionally increasing the number of subpatches for all
patches while keeping fixed the number of grid points per subpatch, as
described in Section~\ref{subsubsec:subpatches}, or 2)~Modifying the
offending patches so as to adequately increase the sizes of the
relevant overlap regions.  An initial assignment of grid sizes aiming
at satisfying the donor-receiver condition might be selected by
ensuring that for each side of each patch $\Omega^{\mathcal{R}}_p$
that is not contained in $\Gamma$, a $2n_v + 1$-point wide layer of
grid points adjacent to that side is also included in a union of one
or more patches $\Omega^{\mathcal{R}'}_{p'}$ (with
$(\mathcal{R}', p') \neq (\mathcal{R}, p)$.

\subsubsection{\label{subsubsec:intrapatch} Intra-patch communication}

\begin{figure}[H]
\centering
    \includegraphics[width=1\linewidth,]{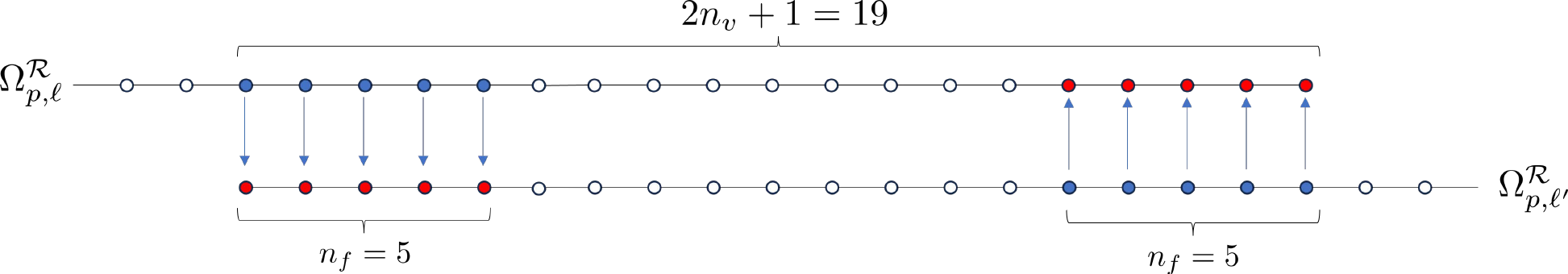}
    \caption{One dimensional illustration of intra-patch communication
      between two subpatches $\Omega^{\mathcal{R}}_{p, \ell}$ and
      $\Omega^{\mathcal{R}}_{p, \ell'}$ of the same patch . Red dots:
      receiver points in the fringe regions of subpatches
      $\Omega^{\mathcal{R}}_{p, \ell}$ and
      $\Omega^{\mathcal{R}}_{p, \ell'}$. Blue dots: donor points in
      subpatch $\Omega^{\mathcal{R}}_{p, \ell}$ (resp.
      $\Omega^{\mathcal{R}}_{p, \ell'}$) used for solution data
      communication with the receiving points of subpatch
      $\Omega^{\mathcal{R}}_{p, \ell'}$ (resp.
      $\Omega^{\mathcal{R}}_{p, \ell}$).}
    \label{fig:intrapatch}
\end{figure}

Unlike the communication between pairs of adjacent patches considered
in Section~\ref{subsubsec:interpatch}, which proceeds via
interpolation from one subpatch fringe region to another, the
``intra-patch'' communication, namely, communication of data among
pairs of neighboring subpatches of a single patch, relies on the set
of fringe points that are {\em shared} between such neighboring
subpatches---in such a way that communication amounts merely to
exchange of solution values as illustrated in
Figure~\ref{fig:intrapatch}. As in the inter-patch case, intra-patch
solution communications are subject to the donor-receiver condition
outlined in Remark~\ref{overlap_cond}. Clearly, the donor-receiver
condition is satisfied for a given patch $\Omega^{\mathcal{R}}_p$
provided every pair $\Omega^{\mathcal{R}}_{p, \ell}$ and
$\Omega^{\mathcal{R}}_{p, \ell'}$
($\ell, \ell' \in \Theta^{\mathcal{R}}_p$) of adjacent subpatches
shares a $2n_f$-point wide layer of grid points along the side common
to both subpatches, where the notion of side of a subpatch is
introduced in Section~\ref{subsec:overlap}. We note that, in view of
the patch decomposition procedures described in Section~\ref{square}
(according to which subpatches of the same patch share a
$(2n_v + 1)$-point wide overlap region with their neighbors), the
subpatch donor-receiver condition is satisfied for the parameter
values $n_f = 5$ and $n_v = 9$ used in this paper.

\section{\label{subsec:time_marching}{FC-based time
    marching}}

\subsection{\label{subsec:differentiation} {Spatio-temporal
    discretization} }

The proposed algorithm spatially discretizes the solution vector
$\mathbf{e}= \textbf{e}(\mathbf{x},t)$ on the basis of $q$-dimensional
families of vector grid functions
\[
  b:\mathcal{G}\to \mathbb{R}^q 
\]
with $q = 4$.  Here, with reference to~\eqref{eq:3_grids}, the overall
spatial computational grid is given by
\begin{equation} 
  \label{eq:overall_grid}
  \mathcal{G} = \bigcup_{\mathcal{R}\in \mathbb{T}}\; \bigcup_{1\leq p\leq P_{\mathcal{R}}}\; \bigcup_{\ell\in \Theta^{\mathcal{R}}_p}
    \mathcal{G}^{\mathcal{R}}_{p, \ell}, \qquad \mathbb{T} =  \{\mathcal{S},\mathcal{C},\mathcal{I} \}.
  \end{equation}
  For each vector grid function $b$, further, we write
\[
b = (b^{\mathcal{R}}_{p, \ell})\quad\mbox{and}\quad  b^{\mathcal{R}}_{p, \ell, i, j} = b^{\mathcal{R}}_{p,
    \ell}(\mathbf{x}_{ij}), \quad\mbox{where}\quad \mathbf{x}_{ij} =
  \mathcal{M}^{\mathcal{R}}_p(q^{\mathcal{R}, 1}_{p, i}, q^{\mathcal{R}, 2}_{p, j}), \quad (i, j) \in \mathcal{D}^{\mathcal{R}}_{p, \ell},
\]
in terms of the parameter space grid points
$q^{\mathcal{R}, 1}_{p, i}$ and $ q^{\mathcal{R}, 2}_{p, j}$ defined
in~\eqref{eq:param-grid}, where $\mathcal{D}^{\mathcal{R}}_{p, \ell}$
is defined in~\eqref{D_R}. In particular, a spatially-discrete but
time-continuous version $\textbf{e}_h = \textbf{e}_h(t)$ of the
solution vector $\textbf{e}(\textbf{x}, t)$ is obtained as a $q=4$
vector grid-function family
$\mathbf{e}_h(t) = (\textbf{e}^{\mathcal{R}}_{h, p, \ell}(t))$ over
all relevant values of $\mathcal{R}$, $p$ and $\ell$; in accordance
with the conventions above we may write
$\textbf{e}^{\mathcal{R}}_{p, \ell,i,j}(t) =
\textbf{e}^{\mathcal{R}}_{h, p, \ell}(\mathbf{x}_{ij},t)$.

Like~\cite[Sec. 3.1]{bruno2022fc}, the proposed algorithm produces the
necessary spatial derivatives on the basis of FC-based
differentiation. In detail, approximations of the partial derivatives
of $\textbf{e}$ with respect to the parameter-space variables are
obtained via the sequential application of one-dimensional
FC-differentiation (Section~\ref{fc}) in the $q^1$ or $q^2$ parameter
variables (see Section~\ref{square}) on the basis of the
point values
$\textbf{e}^{\mathcal{R}}_{p, \ell,i,j}(t) =
\textbf{e}^{\mathcal{R}}_{h, p, \ell}(
\mathcal{M}^{\mathcal{R}}_p(q^{\mathcal{R}, 1}_{p, i}, q^{\mathcal{R},
  2}_{p, j}),t)$. Using such 1D discrete spatial differentiation
operators to obtain each one of the $\textbf{q}$ derivatives in the
curvilinear equation~\eqref{eq:general_curvilinear}, the resulting
discrete equation may be expressed in the form
\begin{equation}\label{semi-discrete}
  \frac{d \textbf{e}_h(t)}{dt} = L[\textbf{e}_h(t)],
\end{equation}
where $L = L^{\mathcal{R}}_{p, \ell}$ denotes a family of discrete
operators that incorporate the products of the Jacobian matrices and
divergences of the viscous and convective fluxes
in~\eqref{eq:general_curvilinear}.

Following~\cite{schwander2021controlling,bruno2022fc}, on the other
hand, the algorithm's time-stepping proceeds via the 4-th order strong
stability preserving Runge-Kutta scheme~\cite{gottlieb2005high}
(SSPRK-4).  This scheme, which leads to low temporal dispersion and
diffusion in smooth flow regions, is employed in conjunction with an
adaptive time-step $\Delta t$ selected at each time $t=t_n$ according
to the expression
\begin{equation} \label{eq: CFL} \Delta t = \frac 1\pi \min_{(\mathcal{R}, p)}
   \min_{\ell \in \Theta^{\mathcal{R}}_p}\frac{\textrm{CFL}}{
    \max_{\textbf{x} \in \mathcal{G}^{\mathcal{R}}_{p,
          \ell}}\left\{\lvert S[\textbf{e}](\mathbf{x}, t)
      \rvert/\widetilde{h}_{\mathcal{R}, p}\right \} +  
      \max_{\textbf{x} \in \mathcal{G}^{\mathcal{R}}_{p, \ell}}\left\{
      \mu[\textbf{e}]\textbf{x},
      t)/(\widetilde{h}_{\mathcal{R}, p})^2\right\}},
\end{equation}
which generalizes the corresponding single patch version used
in~\cite{schwander2021controlling,bruno2022fc}. Here
$\widetilde{h}_{\mathcal{R}, p}$ denotes the minimum mesh size in the
patch $\Omega^{\mathcal{R}}_p$ (a lower-bound close to that minimum
can be used, which can be obtained via consideration of the patch
parametrization gradients), $\textrm{CFL}$ denotes a problem-dependent
constant parameter selected to ensure stability, and
$S[\textbf{e}](\mathbf{x}, t)$ and $\mu[\textbf{e}](\mathbf{x},t)$
denote the {\em maximum wave speed bound} (MWSB) and {\em artificial
  viscosity} operators introduced in
Sections~\ref{subsec:subpatch_visc}
and~\ref{subsec:window_propagation}.

\subsection{Enforcement of boundary conditions and overlap
  communications}\label{bc}
Following the ``conventional method'' described
in~\cite{carpenter1995theoretical}, boundary conditions are enforced
at each intermediate stage of the Runge-Kutta scheme, and, both
Dirichlet and Neumann boundary conditions are incorporated as part of
the differentiation process described in the previous
section. Dirichlet boundary conditions at the time-stage $t_{n,\nu}$
($t_n<t_{n,\nu}\leq t_{n+1}$) for the $\nu$-th SSPRK-4 stage
($\nu = 1, \dots, 5$) of the time-step starting at $t=t_n$, are simply
imposed by overwriting the boundary values of the unknown solution
vector $\textbf{e}_h$ obtained at time $t=t_{n,\nu}$ with the given
boundary values at that time, prior to the evaluation of the spatial
derivatives needed for the subsequent SSPRK-4 stage. Neumann boundary
conditions, in turn, are enforced by constructing appropriate Fourier
continuations after each stage of the SSPRK-4 scheme on the basis of
the modified pre-computed matrix $\widetilde Q$ mentioned in
Section~\ref{fc}; see also Remark~\ref{adiabatic} below. Similarly,
patch/subpatch communication is enforced at every SSPRK-4 stage, by
enacting both the inter- and intra-patch data exchanges described in
Section~\ref{subsec:connectivity} and overwriting the corresponding
solution values. In particular, we have found that imposing
patch/subpatch communication at every stage leads to smoother
numerical solutions, particularly at patch-boundary regions, than are
obtained from a single enforcement at the end of every SSPRK-4 time
step.

\begin{remark}\label{adiabatic_th}
  An important point concerns the assignment of artificial viscosity
  at and around obstacle boundaries, which occurs naturally as part of
  the FC-SDNN algorithm (see e.g. the insets in the right three panels
  in Figure~\ref{Cylinder_Flow_solutions}), and the connection between
  such assignments and the boundary conditions used. (Note from
  Section~\ref{artvisc} that the ``artificial viscosity'' concept used
  throughout this paper commingles a viscosity quantity associated
  with the fluid flow as well as a quantity akin to a heat
  conductivity used in the energy equation.)  Indeed, in order to
  ensure well-posedness in presence of viscosity in a neighborhood of
  the obstacle boundaries, enforcement of no-slip/adiabatic boundary
  conditions (that is, vanishing boundary conditions for both the
  tangential and normal velocity components as well as the normal
  temperature gradient and the associated energy condition introduced
  in Remark~\ref{adiabatic}), is necessary. The FC-SDNN algorithm thus
  produces natively a setting that, as illustrated in
  Section~\ref{subsubsec:flow_cylinder}, incorporates both inviscid
  flow away from boundaries as well as boundary layers near boundaries
  as proposed in Prandtl's boundary-layer
  theory~\cite{prandtl1938berechnung}. We have additionally found
  that, as a significant by-product of this overall scenario for
  treatment of viscosity, boundary layers and boundary
  conditions---and unlike other
  approaches~\cite{nazarov2017numerical,guermond2018second,chaudhuri2011use,mazaheri2019bounded}---the
  present method does not require the use of positivity-preserving
  flux limiters or other artificial devices to prevent the formation
  of negative-density and pressure regions. We hypothesize that such
  unphysical, instability inducing, negative-quantity numerical
  artifacts, for which flux limiters provide a remedy, indicate a
  failure of well-posedness associated with use of incomplete boundary
  conditions in an effectively (numerically) viscous
  environment---e.g. as a constant velocity initial condition
  mismatches vanishing normal velocity conditions on the aft
  obstacle-boundary portions. All of the numerical results presented
  in this paper (Sections~\ref{subsubsec:flow_cylinder}
  through~\ref{subsubsec:shock_cylinder}) have been produced without
  use of positivity preserving limiters.
\end{remark}

\begin{remark}\label{adiabatic}
  With reference to Remark~\ref{adiabatic_th}, the FC-SDNN algorithm
  enforces adiabatic boundary conditions (vanishing of the normal
  gradient of the temperature $\theta$) at all obstacle boundaries (as
  illustrated in Section~\ref{sec:numerical results}). Since the
  temperature is not one of the unknowns utilized in
  equation~\eqref{eq: euler 2d equation}, the enforcement of an
  adiabatic condition requires consideration of
  equation~\eqref{perfect_gas_T}.  The proposed algorithm thus
  proceeds as follows: 1)~At the beginning of each SSPRK-4 time stage,
  the solution $\textbf{e}_h$ is used in conjunction
  with~\eqref{perfect_gas_T} to obtain the grid values of $\theta$ at
  all grid points in the interior of the computational domain, and
  then the adiabatic boundary condition is employed to obtain the
  values of $\theta$ at the physical boundaries by producing, on the
  basis of the method described in~\cite[Sec. 6.3]{amlani2016fc}, a
  Fourier Continuation expansion that matches the interior data and is
  consistent with a vanishing normal gradient of the temperature
  $\theta$; 2)~The Fourier Continuation expansion of $\theta$ obtained
  per point 1 is then used together with the values of the density
  $\rho$ and the horizontal and vertical momenta $\rho u$ and $\rho v$
  that are part of the existing vector $\textbf{e}_h$ (but ignoring
  the values of $E$ that are also part of $\textbf{e}_h$) to
  evaluate the derivatives of $E$ and $p$ with respect to $q^1$ and
  $q^2$ by differentiation of equations~\eqref{p_T}
  and~\eqref{perfect_gas_T} on the basis of the sum and product
  differentiation rules.  Having obtained spatial derivatives of $E$
  and $p$ with respect to $q^1$ and $q^2$ that are consistent with the
  adiabatic boundary conditions, the time stage proceeds by
  3)~Evaluating the spatial derivatives of the remaining components of
  $\textbf{e}_h$ (namely, $\rho$, $\rho u$, $\rho v$) on the basis
  Fourier Continuations of the existing values of these quantities
  under the existing boundary conditions for $u$ and $v$ (no boundary
  conditions are imposed on $\rho$ at physical boundaries where, in
  accordance with the physics of the problem, this quantity is evolved
  in the same manner as at the interior grid points).
\end{remark}

\subsection{\label{subsec:filtering}{Spectral
    filtering}}

Spectral methods regularly utilize filtering strategies in order to
control the error growth in the unresolved high frequency modes, and
the present method is not an exception: like the
approaches~\cite{albin2011spectral,amlani2016fc}, the proposed
algorithm incorporates a multi-patch spectral filter in conjunction
with the Fourier Continuation expansions at the level of each
subpatch. Additionally, a localized initial-data filtering strategy
(which is only applied at the initial time, and then only at
subpatches that contain discontinuities in the initial data), is used
to regularize discontinuous initial conditions, while avoiding
over-smearing of smooth flow-features. (The localized filtering
approach was introduced in~\cite{bruno2022fc} in the context of a
single-grid FC-based shock-dynamics solver.)  Details regarding the
aforementioned multi-patch and localized filtering strategies are
provided in Sections~\ref{global filtering} and~\ref{localized
  filtering}.

\subsubsection{Subpatch-wise filtering strategy \label{global
    filtering}}

To introduce the subpatch-wise filtering method we use we first
consider filtering of a one-dimensional Fourier Continuation expansion
of the form
\[
  \sum_{k = -M}^{M} \hat{F}^c_k \exp (2\pi i k x / \beta), 
\]
for which the corresponding filtered expansion is given by
\begin{equation} \label{eq: filtering}
  \widetilde F = \sum_{k = -M}^{M}
  \hat{F}^c_k \sigma \Big( \frac{2k}{N + C} \Big) \exp (2\pi i k x /
  \beta)
\end{equation}
where
\begin{equation*} \label{eq: filter function} \sigma \Big( \frac{
    2k}{N+C} \Big) = \exp\Big( -\alpha_f \Big(\frac{2k}{N+C}
  \Big)^{p_f} \Big)
\end{equation*}
with adequately chosen values of the positive integer $p_f$ and the
real parameter $\alpha_f>0$. In order to filter a function defined on
a two-dimensional subpatch of the type that underlies the
discretizations considered in this paper, the spectral filter is
applied sequentially, one dimension at a time---thus filtering in
accordance with~\eqref{eq: filtering}, in the $q^1$ subpatch parameter
variable for each value of $q^2$, and then applying the reverse scheme
to the resulting function, filtering with respect to $q^2$ for each
value of $q^1$.

In the strategy proposed in this paper, Fourier Continuation
expansions of the density $\rho$, the horizontal and vertical momenta
$\rho u$ and $\rho v$ and the temperature $\theta$ are computed and
filtered independently on each grid $G^{\mathcal{R}}_{p, \ell}$, at
every time step following the initial time (but not at individual
SSRPK-4 stages), and using the parameter values $\alpha_f = 10$ and
$p_f = 14$. The energy $E$ component of the solution vector
$\textbf{e}_h$ is not directly filtered: per Remark~\ref{adiabatic},
the quantities $E$ and $p$ are re-expressed in terms of the filtered
versions of $\rho$, $\rho u$ and $\rho v$ and $\theta$, and the necessary
derivatives of $E$ and $p$ are obtained via direct differentiation of
equations~\eqref{p_T} and~\eqref{perfect_gas_T}. Thus, the quantities
$E$ and $p$ filtered indirectly, by means of the filters applied to
the other flow variables.  This approach makes it possible to
incorporate the adiabatic boundary conditions (while simultaneously
filtering $E$ and $p$) even in absence of the temperature $\theta$ in the
vector $\textbf{e}$ of unknowns.

\subsubsection{Localized discontinuity-smearing for initial
  data\label{localized filtering}}

To avert the formation of spurious oscillations originating from
discontinuous initial conditions without unduly smearing smooth flow
features, before the initiation of the time-stepping procedure the
initial data is filtered (using a strong form of the spectral
filter~\eqref{eq: filtering}) in all subpatches where initial data
discontinuities exist. This procedure, which generalizes the localized
initial-data filtering approach introduced
in~\cite[Sec. 3.4.2]{bruno2022fc} to the present multi-patch context,
proceeds as follows.

For each subpatch $\Omega^{\mathcal{R}}_{p, \ell}$ within which the
initial condition $\mathbf{e}(\mathbf{x},0)$ is discontinuous, the
localized filtering is performed one dimension at a time in the
parameter-space grid $G^{\mathcal{R}}_{p, \ell}$ and for the physical
quantities $F \in \{\rho, \rho u, \rho v, \theta\}$, in an approach
similar, but different, to the one used in the previous section.  In
detail, considering each restriction of $F$ to the parameter-space
direction $q^i$ ($i=1,2$), containing a discontinuity at a single
point $z$, the smeared-discontinuity function
$\widetilde{F}_\mathrm{sm}$ combines the 1D filtered function
$\widetilde F$ in a neighborhood of the discontinuity with the
unfiltered function elsewhere by means of the expression
\begin{equation}\label{q}
  \widetilde{F}_\mathrm{sm}(x) =  w^{\mathcal{R}, m}_{p, c, r}(x - z) \widetilde F(x) +  (1 - w^{\mathcal{R}, m}_{p, c, r}(x-z)) F(x)\quad, m =1, 2.
\end{equation}
Here, the filtered function $\widetilde F$ is obtained by
applying~\eqref{eq: filtering} with filter parameters $\alpha_f = 10$
and $p_f = 2$, and by using, for $m=1,2$, the window functions
$w^{\mathcal{R}, m}_{p, c, r}$ given by
\begin{equation} \label{eq: qf}
w^{\mathcal{R}, m}_{p, c, r}(x)=
\left\lbrace
  \begin{array}{ccl}
    1 & \mbox{if} &  |x| <  \frac c2h^{\mathcal{R}, m}_p\\
    \cos^2\Big(\frac{\pi (|x| - \frac c2h^{\mathcal{R}, m}_p)}{rh}\Big) & \mbox{if} &    \frac c2 h^{\mathcal{R}, m}_p \leq |x| \leq  (\frac c2 + r)h^{\mathcal{R}, m}_p \\
    0 & \mbox{if} & |x| > (\frac c2 + r)h^{\mathcal{R}, m}_p, \\
  \end{array}\right.
\end{equation}
with $c=18$ and $r=9$.  As in~\cite[Sec. 3.4.2]{bruno2022fc}, in case
multiple discontinuities exist along a given 1D parameter direction in
a given subpatch, the procedure is repeated around each discontinuity
point. Should the support of two or more of the associated windowing
functions overlap, then each group of overlapping windows is replaced
by a single window which equals zero outside the union of the supports
of the windows in the group, and which equals one except in the rise
regions for the leftmost and rightmost window functions in the group.

\section{\label{sec:viscosity}Multi-patch artificial viscosity
  assignment}

The proposed strategy for assignment of artificial viscosity values
$\mu[\textbf{e}]=\mu[\textbf{e}](\mathbf{x},t)$ relies on an extension
of the method described in~\cite[Secs. 3.2, 3.3]{bruno2022fc} to the
multi-patch setting utilized in the present contribution. The proposed
artificial-viscosity algorithm proceeds by producing, at first, a
preliminary subpatch-wise viscosity assignment, on the basis of a
variation of the aforementioned method~\cite{bruno2022fc}; the details
of the new multi-patch algorithm are presented in
Section~\ref{subsec:subpatch_visc}. In order to ensure smoothness of
the artificial viscosity assignment across subpatches (and, in
particular, at the interfaces of neighboring subpatches), a novel
multi-patch windowed-viscosity assignment procedure introduced in
Section~\ref{subsec:window_propagation} is used in conjunction with
the preliminary viscosity assignment mentioned above.

\subsection{\label{subsec:subpatch_visc}Subpatch-wise preliminary
  viscosity assignment}

The subpatch-based viscosity assignment procedure described in what
follows, which determines preliminary artificial viscosity values
$\widehat{\mu}[\textbf{e}](\mathbf{x},t)$ independently for each
subpatch $\Omega^{\mathcal{R}}_{p, \ell}$ (with
$\mathcal{R}=\mathcal{S}$, $\mathcal{C}$, or
$\mathcal{I}$, $1\leq p\leq P_\mathcal{R}$ and
$\ell \in \Theta^{\mathcal{R}}_p$), proceeds on the basis of the
``degree of smoothness'' of a certain ``proxy variable'' function
$\Phi(\textbf{e})(\mathbf{x},t)$ of the unknown solution vector
$\textbf{e}$ on $\Omega^{\mathcal{R}}_{p,
  \ell}$. Following~\cite{bruno2022fc}
and~\cite{schwander2021controlling} we use the proxy variable
$\Phi(\textbf{e})$ equal to the Mach number
$\Phi(\textbf{e}) = \norm {\mathbf{u}} \sqrt{\frac{\rho}{\gamma
    p}}$. As detailed in what follows, utilizing this proxy variable,
a smoothness-classification operator $\tau = \tau[\Phi(\textbf{e})]$
characterizes the degree of smoothness of the function
$\Phi(\textbf{e})$ at a certain time $t$ and over each subpatch
$\Omega^{\mathcal{R}}_{p, \ell}$ by analyzing the oscillations of
restrictions of $\Phi(\textbf{e})$ to certain subsets of the subpatch
grid $\mathcal{G}^{\mathcal{R}}_{p, \ell}$ (where the grid subsets
used are contained in regions including interior portions of the
subpatch $\Omega^{\mathcal{R}}_{p, \ell}$ as well as certain subpatch
regions near physical domain boundaries).

The smoothness-classification operator $\tau$ for each such region is
obtained via consideration of approximations of FC expansions of
$\Phi(\textbf{e})= \Phi(\textbf{e}^{\mathcal{R}}_{p, \ell})$ over
$\Omega^{\mathcal{R}}_{p, \ell}$ (see
Section~\ref{subsec:curvilinear_eqns}) that are obtained from the
discrete numerical solution values
$\bm{\phi}^{\mathcal{R}}_{p, \ell} = \Phi(\mathbf{e}^{\mathcal{R}}_{h,
  p, \ell})$ on $\mathcal{G}^{\mathcal{R}}_{p, \ell}$---and we thus
use the approximation
$\tau[\textbf{e}^{\mathcal{R}}_{p, \ell}] \approx
\tilde\tau[\bm{\phi}^{\mathcal{R}}_{p, \ell}]$ that defines the
discrete operator $\tilde\tau$. For each grid
$\mathcal{G}^{\mathcal{R}}_{p, \ell}$ the discrete operator
$\tilde\tau$ is produced on the basis of an application of the 2D
algorithm presented in~\cite[Sec. 3.2.1]{bruno2022fc}. When applied to
the subpatch $\Omega^{\mathcal{R}}_{p, \ell}$ in the present
multi-patch context, however, the previous
algorithm~\cite[Sec. 3.2.1]{bruno2022fc} is used to obtain smoothness
classification of the solution $\textbf{e}^{\mathcal{R}}_{p, \ell}$
only for a certain ``viscosity-generation subgrid''
$\widetilde{\mathcal{G}}^{\mathcal{R}}_{p,
  \ell}\subset\mathcal{G}^{\mathcal{R}}_{p, \ell}$ described in what
follows.  The fact that smoothness classification is only carried out
over the subgrid $\widetilde{\mathcal{G}}^{\mathcal{R}}_{p, \ell}$
relates to the overlap of patches and subpatches in the discretization
structure we use: in the present context, in order to properly
account for overlaps, the algorithm~\cite[Sec. 3.2.1]{bruno2022fc} is
only applied to stencils of points centered at grid points in the set
$\widetilde{\mathcal{G}}^{\mathcal{R}}_{p, \ell}$. Per the
prescriptions in~\cite[Sec. 3.2.1]{bruno2022fc}, the values of the
operator $\tilde\tau$ are produced numerically on the basis of an
Artificial Neural Network (ANN) whose architecture and training
procedure are described in~\cite[Sec. 3.2.2]{bruno2022fc}. (The ANN
weights and biases are loaded from disc at the FC-SDNN initialization
stage.)

The motivation for the introduction of the viscosity-generation
subgrid $\widetilde{\mathcal{G}}^{\mathcal{R}}_{p, \ell}$ is
twofold. On one hand smoothness classification values, which must be
produced for all patch discretization points
$\mathcal{G}^{\mathcal{R}}_p$ (equation~\eqref{eq:3_grids}), are
produced on the basis of FC expansions at the level of subpatch grids
$\mathcal{G}^{\mathcal{R}}_{p, \ell}$, and, since subpatches overlap,
a decision must be made as to which classification value is used at a
grid point in $\mathcal{G}^{\mathcal{R}}_p$ that belongs to the sets
$\mathcal{G}^{\mathcal{R}}_{p, \ell}$ for two or more values of
$\ell$. The decision is dictated on the basis of accuracy: since the
FC approximation, and, thus, the smoothness classification algorithm,
provide more accurate results at interior points of the subpatch
$\Omega^{\mathcal{R}}_{p, \ell}$ than at points near its boundary, for
$\mathcal{R}=\mathcal{S}$, $\mathcal{C}$ and
$\mathcal{I}$ we define the subgrid
$\widetilde{\mathcal{G}}^{\mathcal{R}}_{p, \ell}$ as the set of all
points in $\mathcal{G}^{\mathcal{R}}_{p, \ell}$ that are not contained
in the fringe region $\mathcal{F}^{\mathcal{R}}_{p, \ell, n_v}$
defined in~\eqref{fringe_region}: 
\begin{equation}\label{viscosity_generation}
  \widetilde{\mathcal{G}}^{\mathcal{R}}_{p, \ell} \coloneqq \mathcal{G}^{\mathcal{R}}_{p, \ell} \setminus \mathcal{F}^{\mathcal{R}}_{p, \ell, n_v}.
\end{equation}
The integer $n_v$ used here equals the one utilized in
Section~\ref{square}; as indicated in that section, the value $n_v=9$
is used throughout this paper. The set of all indices of grid points
in $\widetilde{\mathcal{G}}^{\mathcal{R}}_{p, \ell}$ is denoted by
$\widetilde{\mathcal{D}}^{\mathcal{R}}_{p, \ell}$: using the parameter
space discretization points $q^{\mathcal{R}, 1}_{p, i}$ and
$q^{\mathcal{R}, 2}_{p, j}$ as in equation~\eqref{eq:param-grid}, we
have
\begin{equation}\label{index_viscosity_generation}
  \widetilde{\mathcal{D}}^{\mathcal{R}}_{p, \ell} = \{(i, j) \in \mathcal{D}^{\mathcal{R}}_{p, \ell} \ |\ \mathcal{M}^{\mathcal{R}}_p(q^{\mathcal{R}, 1}_{p, i}, q^{\mathcal{R}, 2}_{p, j})\in  \widetilde{\mathcal{G}}^{\mathcal{R}}_{p, \ell} \}.
\end{equation}

The introduction of the preliminary subpatch-wise viscosity assignment
mentioned in the first paragraph of Section~\ref{sec:viscosity}, which
is the main objective of the present section and is given
in~\eqref{eq: V2D_unscaled}, is based on use of the
viscosity-generation subgrid introduced above together with slight
variants, summarized in what follows, of a number of functions and
operators introduced in~\cite[Sec. 3.3]{bruno2022fc}. In detail, the
definition~\eqref{eq: V2D_unscaled} utilizes
\begin{enumerate}
\item The weight function $R$, which assigns viscosity weights to the
  smoothness classification, and is given by $R(1) = 1.5$, $R(2)=1$,
  $R(3) = 0.5$, and $R(4)=0$;
\item The grid function operator $\widetilde R$ defined by
  ($\widetilde R[\eta]_{ij} = R(\eta_{ij})$);
\item\label{MWSB} The MWSB operator $S[\textbf{e}]$ defined as the upper bound
  $S(\textbf{e}) = |u| + |v| + a$ on the speed of propagation
  $\mathbf{u} \cdot \vec{\kappa} + a$ of the wave corresponding to the
  largest eigenvalue of the 2D Flux-Jacobian (which, in a direction
  supported by the unit vector $\vec{\kappa}$, equals
  $\mathbf{u} \cdot \vec{\kappa} + a$~\cite[Sec. 16.3 and
  16.5]{hirsch1990numerical});
\item The discrete version $\widetilde S[\textbf{e}_h]$ of the
  operator introduced in pt.~\ref{MWSB} above, defined by
\begin{equation}\label{Flux_der_Euler_2}
  \widetilde S[\textbf{e}^{\mathcal{R}}_{h, p, l}]_{ij} = |u^{\mathcal{R}}_{p, i, j}| + |v^{\mathcal{R}}_{p, i, j}|  + a^{\mathcal{R}}_{p, i, j};
\end{equation}
\item The $7\times 7$ localization stencils $L^{i, j}$
  ($(i, j) \in \tilde{\mathcal{D}}^{\mathcal{R}}_{p, \ell}$) (with
  $(\mathcal{R},p, \ell)$-dependence suppressed in the notation)
  defined as the sets of points in
  $\mathcal{G}^{\mathcal{R}}_{p, \ell}$ that surround the grid point
  with index $(i, j)$ in $\mathcal{G}^{\mathcal{R}}_{p, \ell}$. Their
  definition is identical to the localization stencils defined
  in~\cite[Sec. 3.3]{bruno2022fc};
\item The discretization parameter $\widehat{h}_{\mathcal{R}, p}$
  equal to the maximum spacing between two consecutive discretization
  points in the patch $\Omega^{\mathcal{R}}_p$ (an upper-bound close
  to that maximum is used in our implementation, which is obtained via
  consideration of the the distance-amplification factors implicit in
  the patch parametrizations used---such as, e.g., in a circular
  parametrization, the angular distance-amplification factors that
  occur at maximum distance from the circle center).
\end{enumerate}
Using these operators and functions, the preliminary artificial
viscosity operator
$\widehat\mu[\textbf{e}^{\mathcal{R}}_{p, \ell}]$ is defined by
\begin{equation} \label{eq: V2D_unscaled}
  \widehat\mu[\textbf{e}^{\mathcal{R}}_{p, \ell}]_{i,j} =
  \widetilde{R}(\tilde\tau[\bm{\phi}^{\mathcal{R}}_{p, \ell}]_{i,j}) \cdot \max_{(k,
    \ell)\in L^{i,j}}(S[\textbf{e}^{\mathcal{R}}_{p, \ell}]_{k
    \ell}) \widehat{h}_{\mathcal{R}, p}, \quad (i, j) \in \tilde{\mathcal{D}}^{\mathcal{R}}_{p, \ell}.
\end{equation}

\subsection{\label{subsec:window_propagation}Overall viscosity
  operator
  $ \mu^{\mathcal{R}}_{p, \ell} =\mu^{\mathcal{R}}_{p,
    \ell}[\textbf{e}]$}

The overall multi-patch artificial viscosity operator used in this
paper is obtained by exploiting a certain smoothing-blending (SB)
operator $\Lambda$ that smooths the viscosity values provided by the
subpatch-wise preliminary artificial viscosity operator~\eqref{eq:
  V2D_unscaled}, and that, generalizing the single patch
approach~\cite[Sec. 3.3]{bruno2022fc}, additionally combines the
smooth values thus obtained in the various subpatches. The application
of the SB operator ensures a well-defined and spatially smooth
time-dependent artificial viscosity profile over the complete
multi-patch computational domain $\Omega$.

To introduce the SB operator $\Lambda$ we first utilize the window
functions $w^{\mathcal{R}, m}_{p, c, r}$~\eqref{eq: qf} to define, for
$\mathcal{R}=\mathcal{S}$, $\mathcal{C}$ and $\mathcal{I}$, and for
each $(i, j) \in \widetilde{\mathcal{D}}^{\mathcal{R}}_{p, \ell}$
(that is to say, for each $(i, j)$ such that the point
$\mathcal{M}^{\mathcal{R}}_p(q^{\mathcal{R}, 1}_{p, i},
q^{\mathcal{R}, 2}_{p, j})$ lies in the viscosity-generation subgrid
$\widetilde{\mathcal{G}}^{\mathcal{R}}_{p, \ell}$ of
$\mathcal{G}^{\mathcal{R}}_{p, \ell}$), the family of subpatch
windowing functions
\begin{equation} \label{W} W^{\mathcal{R}}_{p, \ell, i, j}(\mathbf{x})
  = w^{\mathcal{R}, 1}_{p, 0, 9}(q^1 - q^{\mathcal{R}, 1}_{p, i})
  w^{\mathcal{R}, 2}_{p, 0, 9}(q^2 -q^{\mathcal{R}, 2}_{p, j}), \quad
  \mathbf{q} = (q^1, q^2) =
  (\mathcal{M}^{\mathcal{R}}_p)^{-1}(\mathbf{x}),\quad \mathbf{x} \in
  \mathcal{G}^{\mathcal{R}}_{p, \ell}.
\end{equation}
We then introduce the family of multi-patch windowing functions
\begin{equation} \label{multipatch_W}
\mathcal{W}^{\mathcal{R}}_{p, \ell, i, j}(\mathbf{x}) =
\left\lbrace
  \begin{array}{ccl}
    W^{\mathcal{R}}_{p, \ell, i, j}(\zeta^{\mathcal{R}}_{p, \ell}(\mathbf{x})) & \mbox{if} &  \mathbf{x} \in \Omega^{\mathcal{R}}_{p, \ell} \\
    0 & \mbox{if} & \mathbf{x} \notin \Omega^{\mathcal{R}}_{p, \ell} \\
  \end{array}\right. ,
\end{equation}
where $\zeta^{\mathcal{R}}_{p, \ell}(\mathbf{x})$, denotes the point
in the real-space grid $ \mathcal{G}^{\mathcal{R}}_{p, \ell}$ that is
``closest'' to $\mathbf{x}$ in the following sense:
$\zeta^{\mathcal{R}}_{p, \ell}(\mathbf{x})$ equals the image of the
point in the grid $G^{\mathcal{R}}_{p, \ell}$ which is the closest to
the pre-image of $\mathbf{x}$ under the patch mapping
$\mathcal{M}^{\mathcal{R}}_p$. In some cases this definition requires
disambiguation, which is achieved as follows: letting
$\mathbf{q} = (q^1, q^2) =
(\mathcal{M}^{\mathcal{R}}_p)^{-1}(\mathbf{x})$ we define
\begin{equation} \label{zeta} \zeta^{\mathcal{R}}_{p,
    \ell}(\mathbf{x}) =
  \mathcal{M}^{\mathcal{R}}_p(\mathrm{rnd} (q^1 / h^{\mathcal{R}, 1}_p) h^{\mathcal{R}, 1}_p, \mathrm{rnd} (q^2 / h^{\mathcal{R}, 2}_p) h^{\mathcal{R}, 2}_p), \quad
  \mathbf{x} \in \Omega^{\mathcal{R}}_{p, \ell},
\end{equation}
where for $t\in\mathbb{R}^+$, $\mathrm{rnd}(t)$ denotes the integer
that is closest to $t$, disambiguated by
$\mathrm{rnd}(n+\frac 12) = n+1$ for every integer $n$.

Using these notations, finally, the normalized multi-patch windowing
functions $\widetilde{\mathcal{W}}^{\mathcal{R}}_{p, \ell, i, j}$, the
multi-patch windowing operator $\Lambda$, and the overall viscosity
operator $\mu^{\mathcal{R}}_{p, \ell}$ are defined by
\begin{equation} \label{multipatch_W_norm}
  \widetilde{\mathcal{W}}^{\mathcal{R}}_{p, \ell, i, j}(\mathbf{x}) =
  \frac{\mathcal{W}^{\mathcal{R}}_{p, \ell, i,
      j}(\mathbf{x})}{\sum_{\mathcal{R}'}\sum_{1 \leq p' \leq
      P_{\mathcal{R}}} \sum_{\ell' \in \Theta^{\mathcal{R}'}_{p'}}
    \sum_{(i', j') \in \widetilde{\mathcal{D}}^{\mathcal{R}'}_{p',
        \ell'}} \mathcal{W}^{\mathcal{R}'}_{p', \ell', i',
      j'}(\mathbf{x})}, \quad \mathbf{x} \in \Omega,
\end{equation}
\begin{equation} \label{Lambda} \Lambda[b](\mathbf{x}) =
  \sum_{\mathcal{R}}\sum_{1 \leq p \leq P_{\mathcal{R}}} \sum_{\ell
    \in \Theta^{\mathcal{R}}_{p}} \sum_{(i, j) \in
    \widetilde{\mathcal{D}}^{\mathcal{R}}_{p, \ell}}
  \widetilde{\mathcal{W}}^{\mathcal{R}}_{p, \ell, i, j}(\mathbf{x})
  b^{\mathcal{R}}_{p, \ell, i, j}\quad\mbox{and}
\end{equation}
\begin{equation} \label{viscosity} \mu^{\mathcal{R}}_{p,
    \ell}[\textbf{e}](\mathbf{x}) =
  \Lambda[\widehat{\mu}(\textbf{e}](\mathbf{x}), \quad \mathbf{x} \in
  \mathcal{G}^{\mathcal{R}}_{p, \ell},
\end{equation}
respectively.

\begin{remark} \label{viscosity_overlap} We note that the viscosity
  operator $\mu^{\mathcal{R}}_{p, \ell}$ is constructed on the basis
  of the preliminary viscosity operator $\widehat{\mu}$, which relies
  on an associated smoothness classification operator $\tilde\tau$---
  that is itself computed over grids
  $\widetilde{\mathcal{G}}^{\mathcal{R}}_{p, \ell}$ contained in
  subpatch interiors, {\em except for subpatches that are adjacent to
    the boundary of $\Omega$}, for which smoothness classification
  near subpatch boundaries cannot be avoided. Numerical experiments
  have shown that favoring smoothness classification values produced
  by stencils located in subpatch interiors helps eliminate spurious
  oscillations that would otherwise arise as shocks or contact
  discontinuities travel from one subpatch to the next. Clearly, the
  use of such interior values requires the existence of sufficiently
  large overlap regions between neighboring subpatches. Moreover, the
  large overlaps allow for the use of window functions with relatively
  wide supports, which in turn facilitates the construction of
  viscosity functions with suitably small gradients.
\end{remark}

\section{\label{sec:algorithm}FC-SDNN algorithm}

Algorithm~\ref{alg_total} presents an outline of the FC-SDNN method
introduced in the previous sections. A fine-grained algorithmic
description of the method, including details concerning multi-core parallelization,
is presented in Part~II.
\begin{algorithm}[H]
  \begin{algorithmic}[1]
    \Begin
    \State \textbackslash \textbackslash Initialization.
    \State Initialize the multi-patch FC-SDNN solver, including 1)~Fetching ANN parameters (Section~\ref{subsec:subpatch_visc}), as well as evaluation of 2)~Multi-patch physical grids $\mathcal{G}^{\mathcal{R}}_{p, \ell}$ (eq~\eqref{eq:3_grids}); 3)~Jacobians $J^{\mathcal{R}, p}_{\mathbf{q}}$ (eq~\eqref{eq:general_curvilinear}); and, 4)~Multi-patch windowing functions $\widetilde{\mathcal{W}}^{\mathcal{R}}_{p, \ell}$ (eq~\eqref{multipatch_W_norm}). 
    \State Initialize time: $t = 0$.
     \State \textbackslash \textbackslash Time stepping.
  \While {$t < T$}  
  \State Evaluate the overall artificial viscosity  $\mu^{\mathcal{R}}_{p, \ell}[\textbf{e}]$ for each subpatch grid (eq~\eqref{viscosity}). 
  \State Filter the solution vector $\textbf{e}^{\mathcal{R}}_{h, p, \ell}$ (Section~\ref{subsec:filtering}).
  \State Evaluate of the temporal step-size $\Delta t$ (eq.~\eqref{eq: CFL}).
  \For {each stage of the SSPRK-4 time stepping scheme}
  \State Evolve the solution vector  $\textbf{e}^{\mathcal{R}}_{h, p, \ell}$ for the SSPRK-4 stage and enforce boundary conditions. (Sections~\ref{subsec:differentiation} and~\ref{bc}.)
  \State\label{sol_comm} Communicate the solution values between neighboring patches and subpatches via exchange and interpolation, as relevant (Sections~\ref{subsubsec:interpatch} and~\ref{subsubsec:intrapatch}).
  \EndFor
  \State Update time: $t = t + \Delta t$
  \State Write solution values to disk at specified time steps $t$.
  \EndWhile
  \End
  \caption{Multi-patch FC-SDNN algorithm}
\label{alg_total}
\end{algorithmic}
\end{algorithm}

\section{\label{sec:numerical results}Numerical results}

This section presents computational results produced by means of a
parallel implementation of the multi-patch FC-SDNN method in a number
of challenging test cases, including problems involving
supersonic/hypersonic flow impinging upon obstacles with smooth
boundaries (cf. Remark~\ref{rem_hyperson}), and including multiple
moving-shocks and contact discontinuities as well as collisions of
strong shocks with obstacles, domain boundaries, contact
discontinuities, other shocks, etc.  It is worthwhile to emphasize
here that, as noted in Remark~\ref{adiabatic_th} and as illustrated in
Section~\ref{subsubsec:flow_cylinder} and
Figure~\ref{Cylinder_Flow_solutions}, and, in fact, by all of the
numerical-examples presented in this paper, but unlike the
methods~\cite{nazarov2017numerical,guermond2018second,chaudhuri2011use,mazaheri2019bounded},
the present approach does not incorporate positivity-preserving flux
limiters or other artificial devices to prevent the formation of
negative-density regions---which, e.g. in the case considered in
Figure~\ref{Cylinder_Flow_solutions}, could manifest themselves in
areas behind the cylinder. The parallel implementation of the
algorithm itself and its extension to problems including non-smooth
obstacles together with a number of tests demonstrating the
algorithm's favorable parallel scaling properties and its ability to
effectively tackle non-smooth obstacle boundaries are presented in
Part~II. The numerical tests presented in the present paper were run
on a Beowulf cluster consisting of 30 dual-socket nodes connected via
HDR Infiniband. Each socket incorporates two 28-core Intel Xeon
Platinum 8273 processors, for a total of 56 cores per node, and 384 GB
of GDDR4 RAM per node. While supported by the Xeon processors, the
hyper-threading capability was not utilized in any of the numerical
examples presented in this paper. The ``production runs'' shown
in Figure~\ref{Cylinder_Flow_solutions} required computing times
ranging from 75 to 250 CPU hours, while those in
Figure~\ref{Cylinder_Shock_solutions} took approximately 13 CPU hours each.
Further details on runtimes and parallel efficiency are presented in
Part~II.


\begin{remark}\label{rem_hyperson}
  The literature on compressible flows at high Mach number, including
  simulation and experiment, is somewhat sparse. An early reference is
  provided by~\cite{billig1967shock}. In particular, this contribution
  indicates that ``real-gas effects will affect the shock-wave shape
  only at Mach numbers 8 and above, and then not too
  significantly''. But the reference also implies that such real-gas
  effects do lead to important deviations on physical observables such
  as density ratios across the shock---including, e.g., a variance
  from an experimental ratio of 14.6 compared to a perfect-gas ratio
  of 5.92 for a flow past a blunt body at Mach number
  19.25. Similarly, the much more recent computational
  reference~\cite{novello2024accelerating} presents numerical
  computations demonstrating significant differences between
  perfect-gas solutions and solutions that incorporate certain
  chemical reactions that are triggered at the high pressures observed
  in high Mach number configurations (a numerical experiment at Mach
  16 on a limited section ahead of the cylindrical obstacle is
  presented in that reference). Perfect-gas flow simulations at high
  Mach numbers have nevertheless been considered repeatedly, including
  the early contribution~\cite{calloway1983pressures} which uses the
  method presented in~\cite{kumar1977numerical} (Navier-Stokes with
  Mach number 22), as well as~\cite{hesthaven2021hybrid} (Mach
  10),~\cite{mazaheri2019bounded} (Mach 10),
  and~\cite{woodward1984numerical} (Mach 10). Thus, perfect gas
  simulations at high Mach numbers remain relevant as mathematical and
  computational testbeds even if they are not physically
  accurate. Possibly, for example, an extension to a multi-species
  context of the FC-SDNN solver proposed in this paper may provide
  similar performance (albeit at the increased cost inherent in the
  multi-species context) as it does in the perfect gas case.
\end{remark}

As indicated above, this section presents numerical results of the
application of the FC-SDNN to a number of physical problems involving
strong shocks as well as supersonic and hypersonic flows. The new
results, here and in Part~II accord with existing numerical results,
theory and experimental data. The results include, in particular,
simulations at very high Mach numbers---for which, nevertheless, good
agreement is observed with certain flow details predicted by previous
theory (shock past a wedge, in Part~II) as well as extrapolations from
experimental data (distance between an obstacle and a reflected bow
shock, in Section~\eqref{subsubsec:flow_cylinder_inf}). Throughout
this section the flows obtained are visualized in two different ways,
namely, via direct contour plots of the density function $\rho$, and
by means of ``Schlieren diagrams''---which,
following~\cite{banks2008study} and~\cite{nazarov2017numerical} we
obtain as plots of the quantity
\begin{equation}\label{exp_schlieren}
  \sigma = \exp \left( -\beta \frac{\lvert \nabla \rho(\mathbf{x}) \rvert - \min_{\mathbf{x} \in \Omega} \lvert \nabla \rho(\mathbf{x}) \rvert }{\max_{\mathbf{x} \in \Omega}\lvert \nabla \rho(\mathbf{x}) \rvert - \min_{\mathbf{x} \in \Omega} \lvert \nabla \rho(\mathbf{x}) \rvert} \right);
\end{equation}
in all the test cases considered the value $\beta = 10$ was used. In
regard to the time-stepping constant $\textrm{CFL}$ in~\eqref{eq:
  CFL} the value $\textrm{CFL} = 0.5$ was used for all the
simulations presented in this paper.


This paper presents numerical illustrations based on two distinct
types of initial conditions. The first type involves the interaction
of a supersonic or hypersonic ``Mach $M$'' flow with an obstacle, with
representative values such as $M = 3$, $3.5$, $6$, $10$, and $25$ (see
Remark~\ref{rem_hyperson}). Such problems are considered in
Sections~\ref{subsubsec:flow_cylinder}
and~\ref{subsubsec:flow_cylinder_inf}; the initial conditions used in
such cases are given by a uniform Mach $M$ flow of the form
\begin{equation}\label{mach_flow_ic}
  (\rho, u, v, p) = (1.4, M, 0 ,1)
\end{equation}
with speed of sound $a = 1$.  The second type of initial conditions
used, which are considered in Section~\ref{subsubsec:shock_cylinder},
corresponds to the interaction of a shock initially located at
$x = x_s$ and traveling toward the region $x \geq x_s$ at a speed $M$
that is supersonic/hypersonic with respect to the speed of sound
$a = 1$ ahead of the shock. For such a Mach $M$ shock the initial
condition is given by
\begin{equation} \label{mach_shock_ic}
(\rho, u, v, p)=
\left\lbrace
    \begin{array}{ccc}
        (\frac{(\gamma + 1) M^2}{(\gamma - 1) M^2 + 2}, \frac{\zeta - 1}{\gamma M}, 0, \zeta) & \mbox{if} & x < x_s   \\ 
        (1.4, 0, 0, 1) & \mbox{if} & x \geq x_s,
    \end{array}\right.
\end{equation}
where $\zeta = \frac{2 \gamma M^2 - \gamma + 1}{\gamma +1}$ denotes
the strength of the shock; see e.g.~\cite{chaudhuri2011use}.

\subsection{\label{subsubsec:flow_cylinder}Supersonic flow past a cylinder in a wind tunnel}

\begin{figure}
  \begin{subfigure}[t]{1\linewidth}
    \centering
    \includegraphics[width=0.525\linewidth]{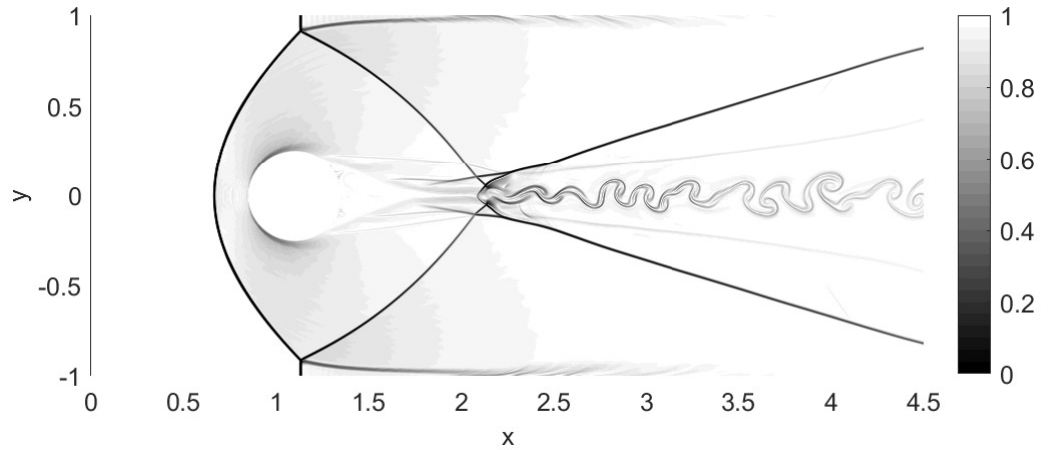}
    \includegraphics[width=0.46\linewidth]{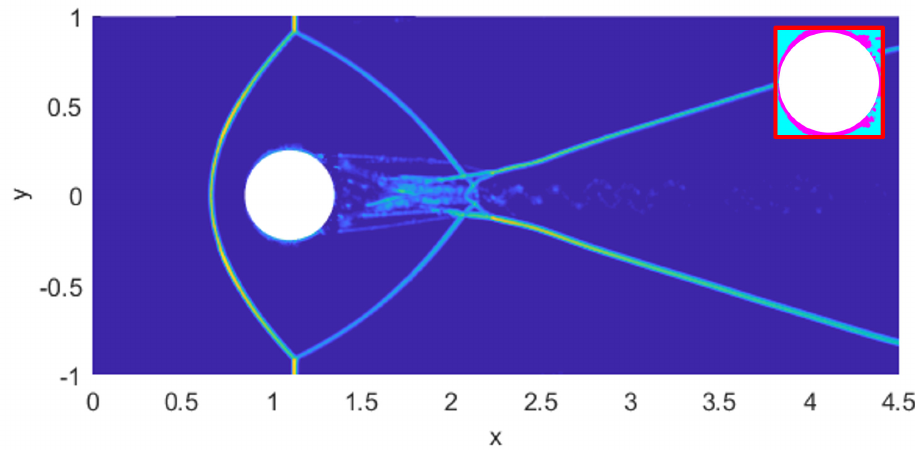}
    \caption{Mach 3 flow.}
  \end{subfigure} 
  \begin{subfigure}[t]{1\linewidth}
    \centering  
    \includegraphics[width=0.525\linewidth]{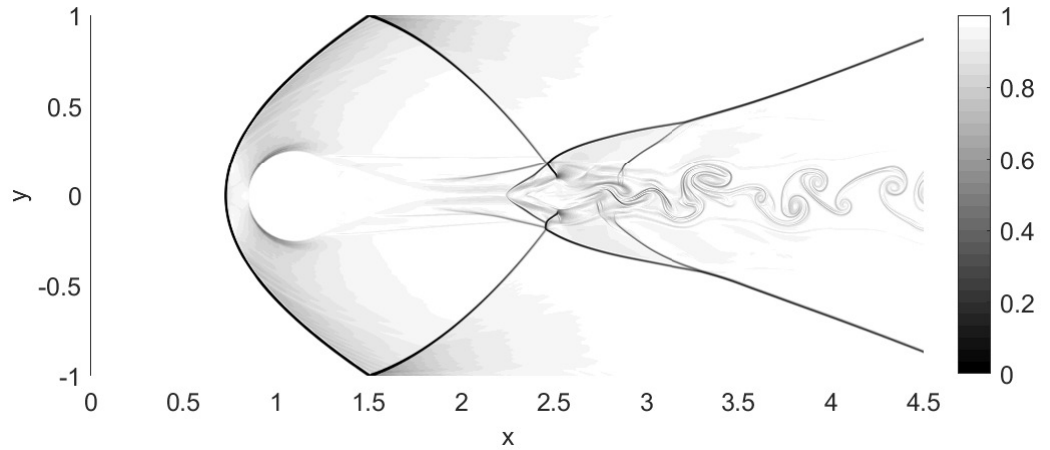}
        \includegraphics[width=0.46\linewidth]{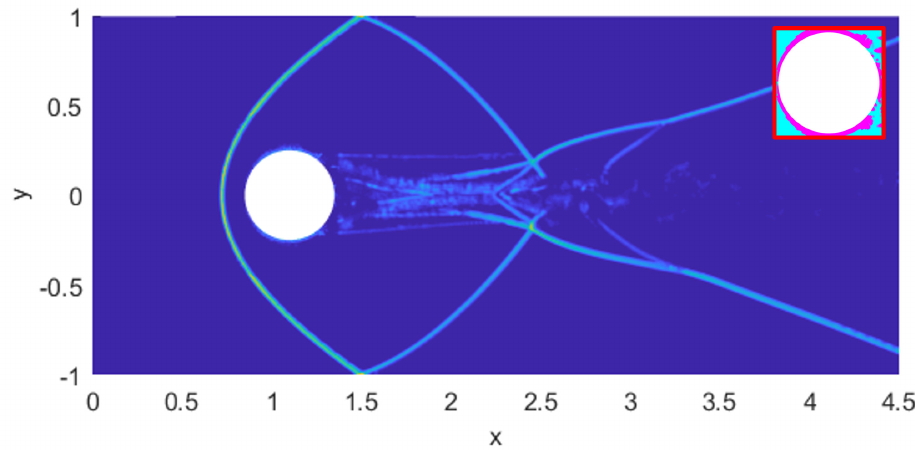}
    \caption{Mach 6 flow.}
  \end{subfigure}

    \begin{subfigure}[t]{1\linewidth}
    \centering  
    \includegraphics[width=0.525\linewidth]{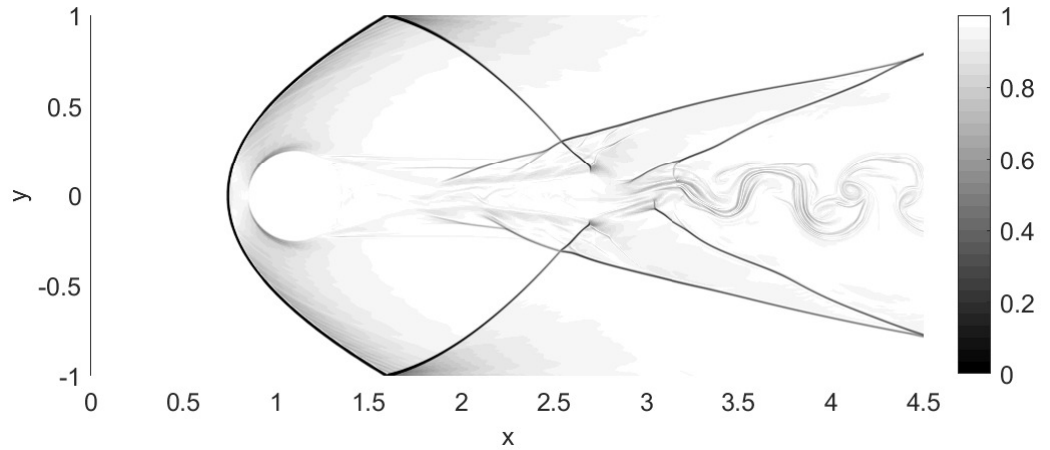}
        \includegraphics[width=0.46\linewidth]{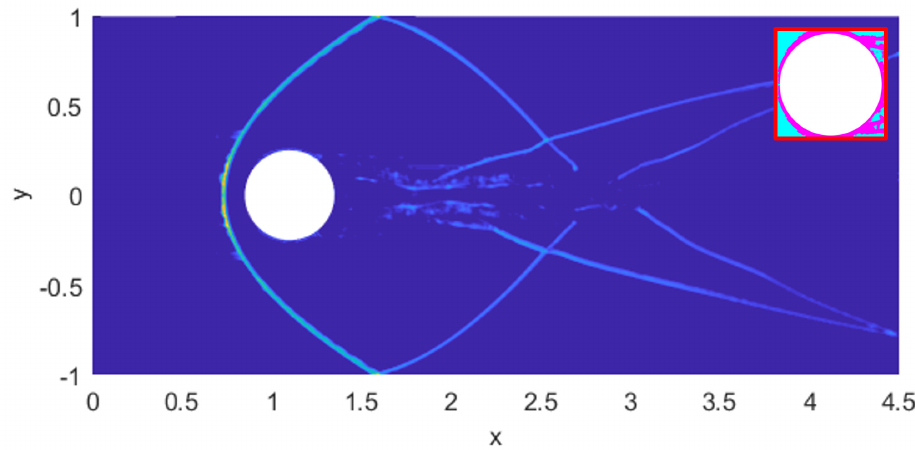}
    \caption{Mach 10 flow.}
  \end{subfigure}

  \caption{Supersonic flows past a cylindrical obstacle at various
    Mach numbers, at time $T = 2$, obtained on a $N = 6.2M$-point
    mesh. Left panels: density Schlieren images. Right panels:
    Artificial viscosity profiles. The insets on the right panels
    display the artificial viscosity profiles (in an appropriate color
    scale) in the immediate vicinity of the cylindrical
    obstacle---which illustrates the boundary-layer resolution
    provided by the method, and which motivates the use of
    no-slip/adiabatic boundary conditions, as discussed in
    Remark~\ref{adiabatic_th}.}
\label{Cylinder_Flow_solutions}
\end{figure}

This section presents results concerning supersonic and hypersonic
flows resulting from initial conditions~\eqref{mach_flow_ic} with Mach
numbers $M = 3$, $M = 6$ and $M = 10$ in a wind tunnel of dimensions
$[0, 4.5]\times[-1, 1]$ over a stationary cylinder of radius
$R_c = 0.25$ and centered at the point $(x_c, y_c) = (1.25, 0)$. An
inflow condition with $(\rho, u, v, p)$ values, coinciding
quantitatively with the initial values, is imposed at the left
boundary at all times. No boundary conditions are imposed on the right
outflow boundary, as befits a supersonic outflow. Reflecting boundary
conditions, corresponding to zero-normal velocity, are imposed at the
bottom ($y = -1$) and top ($y = 1$) walls. No-slip/adiabatic boundary
conditions (Remark~\ref{adiabatic_th}) are imposed at the boundaries
of the cylinder at all times---as befits the viscous-like problem that
is solved in a neighborhood of the cylinder---on account of the
viscosity assigned by the artificial viscosity algorithm, as
illustrated in Figure~\ref{Cylinder_Flow_solutions} and further
discussed in Remark~\ref{adiabatic_th}.

Schlieren diagrams and viscosity profiles of the solutions at $T = 2$
and for various Mach numbers are displayed in
Figure~\ref{Cylinder_Flow_solutions}.  Challenges reported in the
literature for this configuration include difficulties related to the
strong reflected bow shock---characterized by extremely high pressures
between the shock and the cylinder---which frequently leads to spurious
oscillations near the shock. The often observed
formation of vacuum states in the wake region presents further
challenges, for which other approaches employ numerical devices such
as positivity-preserving limiters, particularly for high Mach number
flows, as discussed briefly in the introduction to
Section~\ref{sec:numerical results}. Mach~2, Mach~3, and Mach~3.5
simulations for similar geometric settings were considered
in~\cite{nazarov2017numerical} and~\cite{guermond2018second}
and~\cite{chaudhuri2011use}, respectively.

\subsection{\label{subsubsec:flow_cylinder_inf} $M = 25$ hypersonic
  flow past a cylinder}

\begin{figure}
    \centering
    \includegraphics[width=1\linewidth]{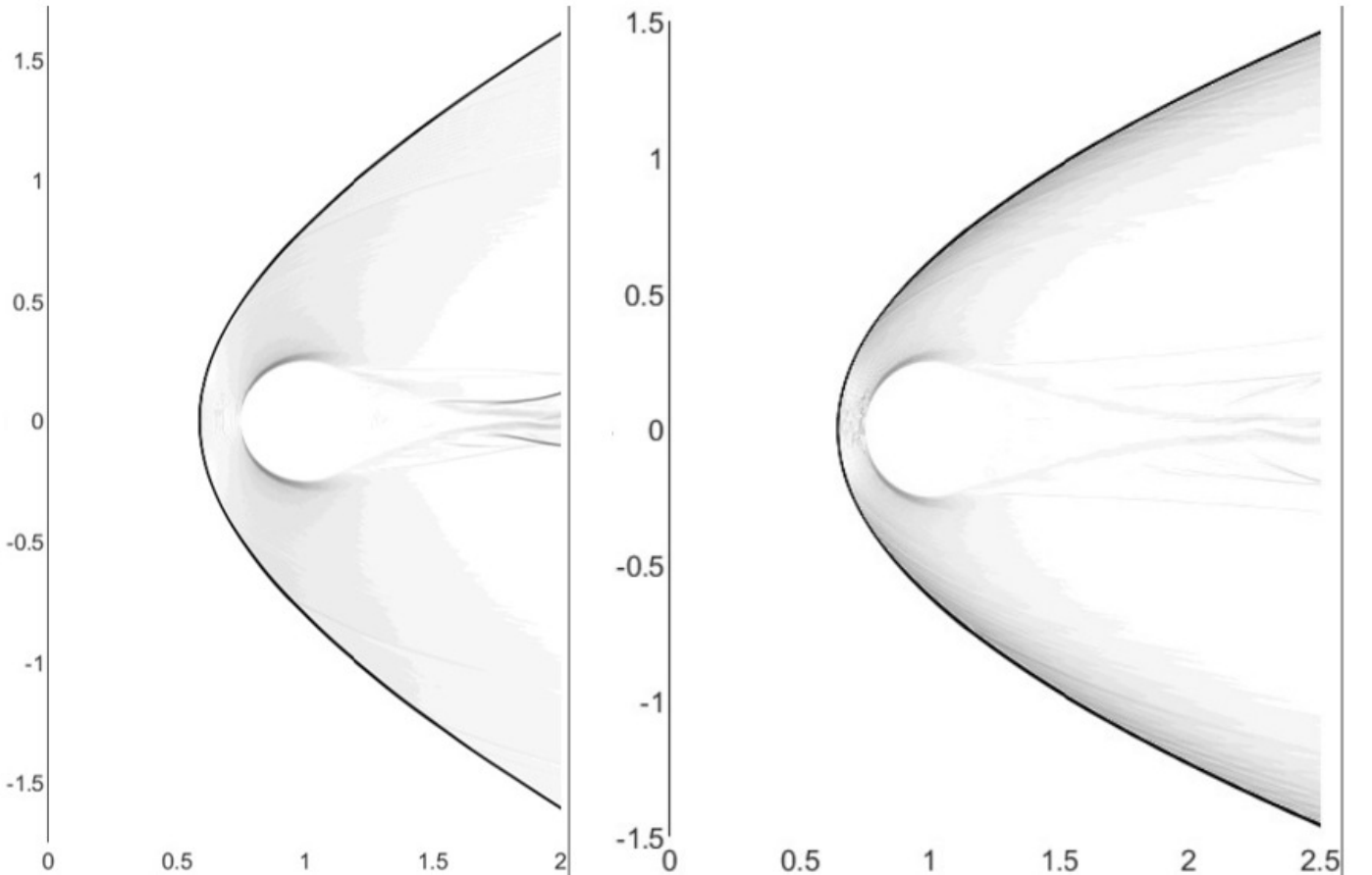}
    \caption{Supersonic flow past a cylindrical obstacle in a wide
      channel, at time $T = 2.0$, for two different Mach numbers, Mach
      3.5 (left) and Mach 25 (right), obtained on meshes containing
      $N \approx 8.8M$ and $N \approx 9.8M$ discretization points,
      respectively.}
\label{Cylinder_wide_Flow_solutions}
\end{figure}

We next consider a set of problems involving flows at supersonic
speeds (Mach 3.5), and re-entry speeds (Mach 25;
cf. Remark~\ref{rem_hyperson}) past a cylindrical obstacle of radius
$R_c = 0.25$ and centered at the point $(x_c, y_c) = (1, 0)$, for
which initial conditions~\eqref{mach_flow_ic} with Mach numbers
$M = 3.5$ and $M = 25$ are imposed, respectively. For the supersonic
(resp. the re-entry) velocity problem the computational domain
considered consists of the portion of the rectangle
$[0, 2]\times[-1.75, 1.75]$ (resp.  $[0, 2.5]\times[-1.5, 1.5]$)
located outside the cylinder. An inflow condition with
$(\rho, u, v, p)$ values coinciding quantitatively with the initial
values is imposed at the left boundary at all times, and no boundary
conditions are imposed on the right boundary, as befits a supersonic
outflow. No-slip/adiabatic (Remark~\ref{adiabatic_th}) boundary
conditions are imposed on the cylinder, while vanishing normal
derivatives for all variables are imposed at the top and bottom domain
boundaries.  Schlieren diagrams of the solutions at $T = 2.0$ are
displayed in Figure~\ref{Cylinder_wide_Flow_solutions}.

A notable feature in these simulations is the formation of a reflected
bow shock ahead of the cylinder. An existing parametric fit to
experiment~\cite{billig1967shock} of the distance $d_0 = d_0(M)$
between the bow shock and the cylinder's leading edge as a function of
the Mach number $M$ is used in what follows (as it was done previously
in~\cite{boiron2009high, chaudhuri2011use}) to validate, at least
partially, the quality of the proposed gas-dynamics solver. According
to~\cite{billig1967shock} the distance $d_0$ is given by
\begin{equation} \label{cylinder_bow} \frac{d_0}{2 R_c} \approx \lambda_1
  \exp \Big(\frac{\lambda_2}{M^2} \Big),
\end{equation}
where, according to~\cite{billig1967shock} we have $\lambda_1 = 0.193$
and $\lambda_2 = 4.67$.  Numerical values of the ratio
$\frac{d_0}{2 R_c}$ produced by our solver for various Mach numbers
are provided in Table~\ref{table:bow}; these data show that, for each
Mach number considered, the numerical value of $d_0$ converges towards
the empirical value given in~\eqref{cylinder_bow} as the mesh is
refined.  The non-asterisked (resp. asterisked) ``Exper. parametric
fit'' quantities in Table~\ref{table:bow} utilize the value
$\lambda_1\approx 0.193= 0.386/2$ provided in the original
reference~\cite{billig1967shock} (resp. the approximated value
$\lambda_1\approx 0.2$ used in~\cite[p. 1738]{chaudhuri2011use}). The
differences between the simulations and the experimental parametric
fit appear consistent with experimental error levels reported
in~\cite{billig1967shock}.

\begin{table}[H]
\centering
\begin{tabular}{|c|c||c|c|}
  \hline
  $M = 3.5$ & $\frac{d_0}{2 R_c}$ &  $M = 25$ & $\frac{d_0}{2 R_c}$  \\ \hline
  Exper. parametric fit  & 0.283 (0.293$^*$)   & Exper. parametric fit  & 0.194 (0.201$^*$)  \\
  $N = 2,203,416$ & 0.32  &   $N = 2,448,240$ & 0.21   \\
  $N = 8,813,664$ & 0.30  &      $N = 9,792,960$  &  0.20    \\
  $N = 19,830,744$ & 0.30  &    $N = 22,034,160$ & 0.20     \\ \hline
\end{tabular}
\caption{Bow-shock/cylinder distance in the supersonic flow past a
  cylinder in a wide channel, including data corresponding to the
  experimental parametric fit~\cite[Fig. 2]{billig1967shock}
  (eq.~\eqref{cylinder_bow} above) and computational results for
  various values of $N$. The non-asterisked (resp. asterisked)
  ``Exper. parametric fit'' quantities in this table utilize the value
  $\lambda_1\approx 0.193= 0.386/2$ provided in the original
  reference~\cite{billig1967shock} (resp. the approximated value
  $\lambda_1\approx 0.2$ used in~\cite[p. 1738]{chaudhuri2011use}).
  The differences between the simulations and the experimental
  parametric fit appear consistent with experimental error levels
  reported in~\cite{billig1967shock}. }
\label{table:bow}
\end{table}

\subsection{\label{subsubsec:shock_cylinder}Shock-cylinder interaction}

This section concerns shock-cylinder interaction problems for the
geometry utilized in Section~\ref{subsubsec:flow_cylinder} with two
different shock speeds of Mach numbers $M = 3$ and $M = 10$. The
corresponding initial conditions are given by~\eqref{mach_shock_ic},
with initial shock position at $x_s = 0.5$.  An inflow condition with
$(\rho, u, v, p)$ values coinciding with the $x\leq x_s$ initial
values is imposed at the left boundary at all times, and an outflow
condition consisting of the time-independent pressure value $p = 1$ is
imposed at the right boundary, also at all times.  Slip-wall boundary
conditions are imposed at the bottom and top walls, while no-slip and
adiabatic boundary conditions (Remark~\ref{adiabatic_th}), finally,
are imposed at the boundaries of the cylinders at all times.

\begin{figure}
  \begin{subfigure}[t]{0.5\linewidth}
    \centering
    \includegraphics[width=1\linewidth]{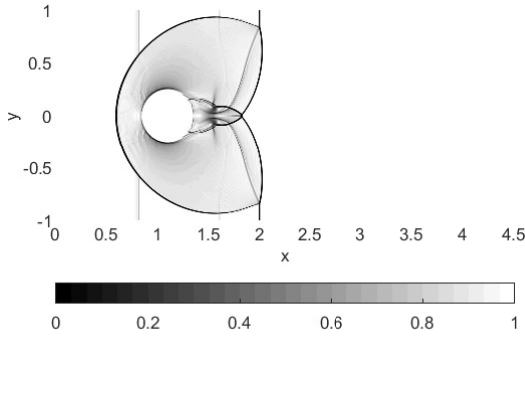}
  \end{subfigure} 
  \begin{subfigure}[t]{0.5\linewidth}
    \centering  
    \includegraphics[width=1\linewidth]{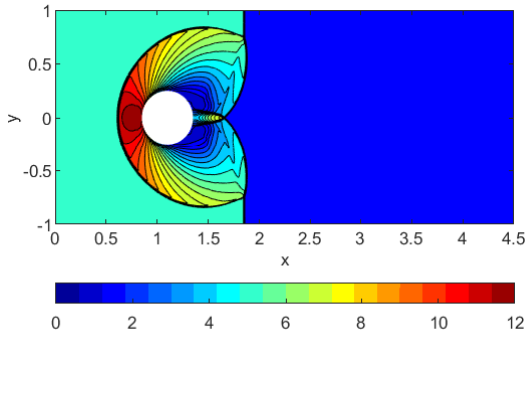}
  \end{subfigure}

  \begin{subfigure}[t]{0.5\linewidth}
    \centering
    \includegraphics[width=1\linewidth]{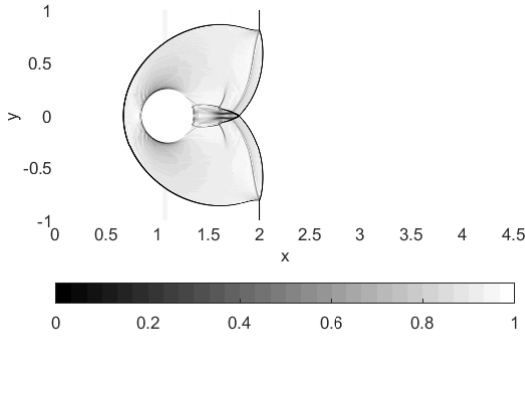}
  \end{subfigure} 
  \begin{subfigure}[t]{0.5\linewidth}
    \centering  
    \includegraphics[width=1\linewidth]{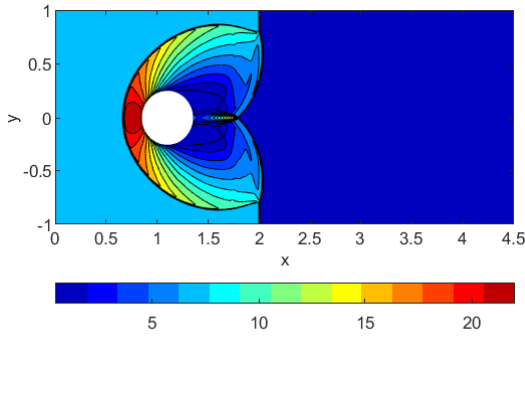}
  \end{subfigure}
  \caption{Shock-cylinder interaction for two different shock speeds
    and at two different times, namely, Mach 3.5 at time $T = 0.45$
    (top) and Mach 10 at time $T = 0.15$ (bottom), obtained on an
    $N \approx 6.2M$-point mesh. Left panels: density Schlieren
    images. Right panels: density contour plots.}
\label{Cylinder_Shock_solutions}
\end{figure}

Schlieren and density contour plots for the Mach 3 solution (resp. the
Mach 10 solution) at time $T = 0.45$ (resp.  $T = 0.15$) are displayed
in Figure~\ref{Cylinder_Shock_solutions}. The Mach 3 Schlieren
diagrams display flow features in close agreement with the experimental
and numerical results presented
in~\cite[Fig. 4]{bryson1961diffraction}
and~\cite[Fig. 17]{chaudhuri2011use}, respectively, including the
reflected shock and symmetric Mach shocks and contact discontinuities
starting at the back of the cylinder and intersecting the incident
shock at two triple points. Similar features are observed for the Mach
10 problem that is also illustrated in
Figure~\ref{Cylinder_Shock_solutions}. (A single shock-cylinder
interaction problem, with Mach 2.8 shock speed, is presented
in~\cite{chaudhuri2011use}.) The contour plots in
Figure~\ref{Cylinder_Shock_solutions} display smooth density contours,
without spurious oscillations.

\section{\label{sec:conclusion}Conclusions}

This paper, Part-I of a two-part contribution, introduced a novel
multi-patch computational algorithm for general gas-dynamics problems,
namely, the FC-SDNN {\em spectral} shock-dynamics solver, which,
without the taxing CFL constraints inherent in other spectral schemes,
is applicable to {\em general domains}, and up to and including {\em
  supersonic} and {\em hypersonic} regimes.  The proposed algorithm
thus vastly enlarges the domain of applicability beyond that of the
previous single-patch/simple-domain FC-SDNN
solver~\cite{bruno2022fc}---enabling accurate solutions, with smoothly
contoured flow profiles away from shocks, even for fast, highly
challenging fluid flow problems. Notably, this paper includes
illustrations at higher Mach numbers than previously reported for
certain important experimental configurations and physical obstacles
(cf. Remark~\ref{rem_hyperson}).  In particular, the new multi-patch
approach 1)~Enables the application of the FC-SDNN method for
evaluation of flows and shock-dynamics in general 2D domains with {\em
  smooth} boundaries, and in the context of both {\em supersonic and
  hypersonic} flows; 2)~It introduces a {\em smooth and localized}
artificial-viscosity quantity of the type utilized in a previously
existing single-domain FC-SDNN strategy, but which is additionally
compatible with the multi-patch general-domain discretization strategy
used presently; 3)~It does {\em not} require the use of
positivity-preserving flux limiters (Remark~\ref{adiabatic_th}); 4)~It
provides a setting that can be leveraged to incorporate an MPI-based
parallel implementation (described in~\cite{bruno_leibo_partII}); and
5)~It exhibits numerical results (in Part~I and Part~II) in close
accordance with {\em physical theory} and prior {\em experimental and
  computational} results, up to and including both the {\em
  supersonic} and {\em hypersonic} regimes.  The following Part~II
extends the proposed method in several ways, enabling applications to
non-smooth domains, and including further validation and benchmarking
as well as details concerning parallel implementation;
extensions including 3D contexts and GPU parallelization are
envisioned for future work.
\section*{Acknowledgments}
The authors gratefully acknowledge support from NSF and AFOSR under
contracts DMS-2109831, FA9550-21-1-0373 and FA9550-25-1-0015.  \small

\bibliographystyle{abbrv}
\bibliography{bibliography}

\end{document}